\documentclass[11pt]{article}
\usepackage{graphicx}
\usepackage[utf8x]{inputenc}
\usepackage[english]{babel}
\usepackage{amsmath}
\usepackage{amsfonts}
\usepackage{amssymb}
\usepackage{amsthm}
\usepackage[margin=2.5cm]{geometry}
\newtheorem{thm}{Theorem}
\newtheorem{proposition}[thm]{Proposition}
\newtheorem{lemma}[thm]{Lemma}
\newtheorem{corollary}[thm]{Corollary}
\newtheorem{remark}[thm]{Remark}
\newtheorem{conj}[thm]{Conjecture}

\newcommand{\R}{\mathbb{R}}
\newcommand{\N}{\mathbb{N}}
\newcommand{\f}{\varphi}
\newcommand{\e}{\varepsilon}
\renewcommand{\L}{\mathcal{L}}
\newcommand{\Epi}{\mathrm{Epi}}
\newcommand{\conv}{\mathrm{Conv}}
\newcommand{\dom}{\mathrm{dom}}

\title{The functional form of Mahler's conjecture for even log-concave functions in dimension $2$
.}
\author{Matthieu Fradelizi, Elie Nakhle}

\begin{document}
\maketitle


\begin{abstract}
Let $\f:\R^{n}\rightarrow\R\cup\{+\infty\}$ be an even convex function and $\L{\f}$ be its Legendre transform. We prove the functional form of Mahler's conjecture concerning the functional volume product $P(\f)=\int e^{-\f}\int e^{-\L\f}$ in dimension 2: we give the sharp lower bound of this quantity and characterize the equality case.
The proof uses the computation of the derivative in $t$ of $P(t\f)$ and ideas due to Meyer \cite{M} for unconditional convex bodies, adapted to the functional case by Fradelizi-Meyer \cite{FM08a} and extended for symmetric convex bodies in dimension 3 by Iriyeh-Shibata \cite{IS} (see also \cite{FHMRZ}). 
\end{abstract}

\section{Introduction}
In the theory of convex bodies, many geometric inequalities can be generalized to functional
inequalities. This is the case of the Prékopa-Leindler inequality, which is the functional form of the Brunn-Minkowski inequality.  
Let us mention also the Blaschke-Santal{\'o} inequality \cite{San}, which states in the symmetric case that if $K$ is a symmetric convex body in $\R^n$ (in this paper, $K$ symmetric means $K=-K$) and 
$$
P(K)=|K||K^*|,
$$
where $K^*=\{y\in\R^{n}; \langle y, x \rangle\leq 1, \mbox{ for all } x\in K\}$ is the polar body of K then
$$
P(K)\leq  P(B_{2}^{n}),
$$
with equality if and only if $K$ is an ellipsoid, where $B^{n}_{2}$ is the Euclidean ball associated to the standard scalar product in $\R^n$ and $|B|$ stands for the Lebesgue measure of a Borel subset $B$ of $\R^n$ (\cite{P}, see \cite{MP} or also \cite{MR} for a simple proof of both the inequality and the case of
equality).\\

Mahler conjectured an inverse form of the Blaschke-Santal{\'o} inequality for symmetric convex bodies in $\R^n$ in \cite{Mah1}. He asked if for every symmetric convex body $K$,
$$
P(K)\geq P([-1,1]^n)=\frac{4^{n}}{n!}.
$$
It was later conjectured that the equality case occurs if and only if $K$ is a Hanner polytope (see \cite{RZ} for the definition).
The inequality was proved by Mahler for $n=2$ \cite{Mah1} (see also \cite{Re,M} and \cite{S14} Section 10.7, for other proofs and the characterization of the equality case). This conjecture has been proved also in a number of particular cases, for zonoids by Reisner \cite{Re} (see also \cite{GMR}) and for unconditional convex bodies by Saint Raymond \cite{SR} (see also \cite{M}), for hyperplane sections of $B_p^n=\{x\in\R^n; \sum |x_i|^p\le1\}$ and Hanner polytopes by Karasev \cite{K}. The $3$-dimensional case of the conjecture was proved by Iriyeh and Shibata \cite{IS} (see \cite{FHMRZ} for a shorter proof).

Some functional versions of the previous inequalities were proposed with convex bodies replaced with log-concave functions, and polarity with the Legendre transform. Let $\f:\R^{n}\to\R\cup\{+\infty\}$ be an even convex function. Then, the Legendre transform $\L\f$ of $\f$ is defined, for $y\in\R^n$, by 
$$
\L\f(y)=\sup_{x\in\R^{n}}(\langle x,y\rangle-\f(x)).
$$
We define the functional volume product of an even convex function to be
$$
P(\f) =\int e^{-\f(x)}dx\int e^{-\L\f(y)}dy.
$$
The functional version of the Blaschke-Santal{\'o} inequality for even convex functions states that
$$
P(\f)\leq P\left(\frac{|\cdot|^{2}}{2}\right)=(2\pi)^n,
$$
where $|.|$ stands here for the Euclidean norm in $\R^{n}$, with equality if and only if $\f$ is a positive quadratic form. This statement was proved by Ball \cite{B} (see also Artstein-Klartag-Milman \cite{AKM}, Fradelizi-Meyer \cite{FM07} and Lehec \cite{L1, L2} for more general results). For $\f(x)=\|x\|_K^2/2$, one has $\L\f(y)=\|y\|_{K^*}^2/2$, thus it is not difficult to see that
\[
P(\f)=\frac{(2\pi)^n}{|B_2^n|^2}P(K).
\]
This shows that the functional form indeed implies the geometric form of the inequality. 
We deal, in this article, with the following functional version of Mahler's conjecture for even convex functions, that was stated in \cite{FM08a}. We prove it for $n=2$.
\begin{conj}\label{conj}
Let $n\ge1$ and $\f:\R^n\to\R\cup\{+\infty\}$ be an even convex function such that $0<\int_{\R^{n}}e^{-\f(x)}dx<+\infty$. Then,
\[
\int e^{-\f}\int e^{-\L\f}\ge 4^n,
\]
with equality if and only if there exists  $c\in\R$ and two Hanner polytopes  $K_1\subset F_1$ and $K_2\subset F_2$, where $F_1$ and $F_2$ are two complementary subspaces in $\R^n$, such that for all $(x_1,x_{2})\in F_1\times F_{2}$
\[\f(x_1+x_2)=c+\|x_1\|_{K_1}+I_{K_2}(x_2),\quad a.e.\]
where $I_{K}$ is the function defined by $I_{K}(x)=0$ if $x\in K$ and $I_{K}(x)=+\infty$ if $x\notin K$.
\end{conj}
For unconditional functions (and in particular if $n=1$), the inequality in conjecture 1 was proved in \cite{FM08a,FM08b} and the equality case was proved in $\cite{FGMR}$. 
For a symmetric convex body $K$ of $\R^{n}$
and for any $y\in\R^{n}$, one has 
$$
\L I_{K}(y)=\sup_{x\in\R^{n}}(\langle x,y\rangle-I_{K}(x))=\sup_{x\in K}\langle x,y\rangle=h_{K}(y)=\|y\|_{K^{*}},
$$
where $h_{K}$ denotes the support function of $K$. In addition, using Fubini, we have
$$
\int_{\R^{n}}e^{-\|y\|_{K^{*}}}dy=\int_{\R^{n}}\int_{\|y\|_{K^{*}}}^{+\infty}e^{-t}dtdy=\int_{0}^{+\infty}\int_{\{\|y\|_{K^{*}}\leq t\}}e^{-t}dydt=\int_{0}^{+\infty}|tK^{*}|e^{-t}dt=n!|K^{*}|.
$$
Thus, if $\f=I_{K}$ we get 
$$
P(\f)=\int_{\R^{n}} e^{-I_{K}}\int_{\R^{n}} e^{-\|y\|_{K^{*}}}=|K|\int_{\R^{n}} e^{-\|y\|_{K^{*}}}=n!|K||K^{*}|=n!P(K).
$$
Hence, Conjecture 1 implies Mahler's conjecture for symmetric convex bodies. Notice that, recently, Gozlan \cite{G} established precise relationships between the functional form of Mahler's conjecture and the deficit in the Gaussian log-Sobolev inequality. Thus, our results imply better bounds in dimension 2 for these deficits.\\
In the proof, as in \cite{IS, FHMRZ}, we use the notion of equipartition. Denote by $(e_{1},...,e_{n})$ the canonical basis of $\R^{n}$. A function $\f:\R^{n}\to\R\cup\{+\infty\}$ is equipartitioned  if 
$$
\int_{\R^{n}_{\varepsilon}}\f e^{-\f}=\frac{1}{2^{n}}\int_{\R^{n}}\f e^{-\f}\quad\mbox{and}\quad\int_{\R_{\varepsilon}^{n}}e^{-\f}=\frac{1}{2^{n}}\int_{\R^{n}}e^{-\f},
$$
where $\forall\varepsilon\in\{-1,1\}^{n}$, $\R^{n}_{\varepsilon}=\{x\in\R^{n}; \varepsilon_{i}x_{i}\geq 0,\forall i\in\{1,...,n\}\}$.
We use the fact that, in dimension $n\le2$, every even convex function  has an equipartitioned position, see Lemma~\ref{lem:equip}. Moreover, we also prove, in our Theorem \ref{thm:dimn}, that if one has an even convex function $\f$ on $\R^n$ such that $\f$ and $\f_i=\f_{|e_i^\bot}$ are equipartitioned for all $1\le i\le n$ and $\f_i$ satisfy the inequality of the conjecture in dimension $n-1$ then $\f$ satisfies the conjectured inequality in dimension $n$.

This paper is organized in the following way. In section 2, we present some general results on the Legendre transform. In section 3, we establish some properties of the functional volume product. In section 4, we apply the results of sections 2 and 3 to prove the inequality and the case of equality of the functional volume product of $\f$ in dimension 2. Finally, in section 5, we prove the inequality in  dimension $n$ for "strongly" equipartitioned convex functions.

\section{General results on convex functions and  Legendre transform}

Let us recall some useful facts about convex functions and the Legendre duality that can be found in parts $\rm{I}$ and $\rm{V}$ of the book of Rockafellar \cite{R}.
Recall that if $K$ is convex, closed and contains $0$ then $(K^{*})^{*}=K$.  Let $\f :\R^n\to\R\cup\{+\infty\}$ be a convex function. We denote the domain of $\f$ by
\[
\dom(\f)=\{x\in\R^n; \f(x)<+\infty\}.
\]
It is a convex set. If $\f$ is moreover lower semi-continuous and $\dom(\f)\neq\emptyset$ then $\L\L\f=\f$.
The following lemma recalls some standard facts that can be found for example in Lemma 4 of \cite{G}.
\begin{lemma}\label{lem:integ}
Let $n\ge1$ and $\f:\R^n\to\R\cup\{+\infty\}$ be a convex function such that $\min\f=\f(0)=0$. 
Then, the following are equivalent.
\begin{enumerate}
    \item One has $0<\int_{\R^{n}}e^{-\f(x)}dx<+\infty$.
    \item The set $K_\f:=\{x\in\R^{n}, \f(x)\leq 1\}$ is convex bounded and contains $0$ in its interior.
    \item There exists $a,b>0$ such that for every $x\in\R^n$ one has $a|x|-1\le\f(x)\le I_{bB_2^n}(x)+1$.
\end{enumerate}
\end{lemma}
Notice that for every $x\in\R^n$ one has $a|x|-1\le\f(x)\le I_{bB_2^n}(x)+1$ if and only if for every $y\in\R^n$ one has $b|y|-1\le\L\f(y)\le I_{aB_2^n}(y)+1$. Thus, $0<\int_{\R^{n}}e^{-\f}<+\infty$ is equivalent to $0<\int_{\R^{n}}e^{-\L\f}<+\infty$.\\

We define the analogue of sections and projections of convex sets for convex functions. The section of $\f$ by an affine subspace $F$ is simply the restriction of $\f$ to this subspace and is denoted by $\f_{|F}$. The projection $P_F\f: F\to \R\cup\{+\infty\}$ of $\f$ onto a linear subspace $F$ is defined, for $x\in F$, by
\[
P_F\f (x)=\inf_{z\in F^\bot} \f(x+z).
\]
The term projection comes from the fact that if  $\tilde{P}_F:\R^n\times\R\to F\times\R$ denotes the orthogonal projection on $F\times\R$ parallel to $F^\bot$ then $\tilde{P}_F(\Epi(\f))=\Epi(P_F\f)$, where $\Epi(\f)=\{(x,t)\in\R^{n}\times\R;\f(x)\leq t\}$. The infimal convolution of two convex functions $\f$,$\psi:\R^{n}\to\R\cup\{+\infty\}$ is
$$
\f\square\psi(x)=\inf_{z\in\R^n}(\f(x-z)+\psi(z)).
$$
The infimal convolution is a convex function
and it interacts with the Legendre transform in the following way: $\L(\f\square\psi)=\L\f+\L\psi$.
One can also define the projection using infimal convolution by noticing that, for any $x\in\R^n$,
\[
\f\square I_{F^\bot}(x)=\inf_{z\in\R^n} (\f(x-z)+I_{F^\bot}(z))=\inf_{z\in F^\bot} \f(x-z).
\]
Hence, for $x\in F$, one has $\f\square I_{F^\bot}(x)=P_F\f (x)$. Thus, the same nice duality relationship between sections and projections, that holds for convex sets, holds also for convex functions: for any $y\in F$ 
\begin{eqnarray*}
P_F(\L\f)(y)&=&\L\f\square I_{F^\bot}(y)=\L\f\square\L I_F(y)=\L(\f+I_F)(y)=\sup_{x\in\R^n}(\langle x,y\rangle-\f(x)-I_F(x))\\
&=&\sup_{x\in F}(\langle x,y\rangle-\f(x))=\L(\f_{|F})(y),
\end{eqnarray*}
where, by an abuse of notation, we have denoted in the same way by $\L$ the Legendre transform applied to a function defined on $\R^n$ or on a subspace $F$. In each situation, the supremum in the Legendre transform should be understood as taken in the subspace where the function is defined.

For $1\le i\le n$, we denote the restriction of $\f$ to $e_i^\bot$ by $\f_i=\f_{|e_i^\bot}$ and we define the analogue of the projection onto $e_i^\bot$ to be the function $P_i\f:e_i^\bot\to\R\cup\{+\infty\}$ defined, for $x\in e_i^\bot$, by 
\[ P_i\f(x)=P_{e_i^\bot}\f(x)=\inf_{t\in\R}\f(x+te_i).\]
From the preceding, for every $y\in e_i^\bot$, one has 
\[
\L\f_i(y)=\sup_{x\in e_i^\bot}(\langle x,y\rangle -\f(x))=P_i\L\f(y)=\inf_{t\in\R}\L\f (y+te_i). 
\]
\begin{lemma}\label{Param}
Let $\f:\R^{n}\to\R\cup\{+\infty\}$ be differentiable and strictly convex such that $0<\int e^{-\f}<+\infty$, and for $1\leq i\leq n$ let $\f_{i}=\f_{|e^{\perp}_{i}}$ then 
\begin{enumerate}
    \item The function $\nabla\f$ is a bijection from $\R^{n}$ to $\dom(\L\f)$.
    \item One has $\nabla\f(e_{i}^{\perp})=\{y+t_i(y)e_i; y\in\dom(\L\f_{i})\}$ where $t_{i}(y)=\langle\nabla\f\circ(\nabla\f_{i})^{-1}(y),e_{i}\rangle$.
\end{enumerate}
\end{lemma}
\begin{proof}
Since $\f$ is differentiable, we deduce from \cite{R} Theorem 26.3 that $\L\f$ is strictly convex.\\
1. The fact that $\nabla\f$ is a bijection can be found in Corollary 26.3.1 in \cite{R}.\\
2. Since the supremum $\L\f(y)=\sup(\langle x,y\rangle-\f(x))$ is reached at $x=(\nabla\f)^{-1}(y)$ one has 
\[
\L\f(\nabla\f(x))=\langle x,\nabla\f(x)\rangle -\f(x)
\]
and one can conclude from Corollary 23.5.1 of \cite{R} that $(\nabla\f)^{-1}=\nabla(\L\f)$.
Let now $y\in\dom(\L\f_{i})$ be fixed and $g_y(t)=\L\f(y+te_i)$. The function $g_y$ is strictly convex and tends to infinity at infinity so there exists a unique $t_i(y)\in\R$ at which the function $g_y$ reaches its infimum and it satisfies $g_y'(t_i(y))=0$, {\em i.e.} $\langle\nabla\L\f(y+t_i(y)e_i),e_i\rangle=0$ which means that $(\nabla\f)^{-1}(y+t_i(y)e_i)=\nabla\L\f(y+t_i(y)e_i)\in e_i^\bot$. This also means equivalently  that $t_i(y)$ is the unique $t\in\R$ such that $y+te_i\in\nabla\f(e_i^\bot)$. Hence, $\nabla\f(e_i^\bot)=\{y+t_i(y)e_i; y\in\dom(\L\f_{i})\}$ and the orthogonal projection $P_i$ onto $e_i^\bot$ is a bijection from $\nabla\f(e_i^\bot)$ onto $e_i^\bot$. 
Moreover, one has 
\[
P_i\L\f(y)=\inf_{t\in\R}\L\f (y+te_i)=\inf_{t\in\R} g_y(t)=g_y(t_i(y))=\L\f(y+t_i(y)e_i).
\]
Thus,
\[
P_i\L\f(y)=\sup_{x\in\R^n}(\langle x,y+t_i(y)e_i\rangle-\f(x))\ge\sup_{x\in e_i^\bot}(\langle x,y+t_i(y)e_i\rangle-\f(x))=\sup_{x\in e_i^\bot}(\langle x,y\rangle-\f(x))=\L(\f_i)(y).
\]
But in fact, we know that, in the above, the left hand side supremum is reached at \[x=(\nabla\f)^{-1}(y+t_i(y)e_i)=\nabla\L\f(y+t_i(y)e_i)\in e_i^\bot,
\]
hence the above inequality is an equality. 
But the right hand side supremum is reached at $x=(\nabla\f_i)^{-1}(y)$. Since they are reached at the same point this implies that $(\nabla\f)^{-1}(y+t_i(y)e_i)=(\nabla\f_i)^{-1}(y)$,  and $y+t_i(y)e_i=\nabla\f\circ(\nabla\f_i)^{-1}(y)$, thus 
$t_i(y)=\langle \nabla\f\circ(\nabla\f_i)^{-1}(y),e_i\rangle,$ where $(\nabla\f_{i})^{-1}$ is a bijection from $\dom(\L\f_{i})=\dom P_{i}(\L\f)$ to $e_{i}^{\perp}$. 
\end{proof}

For every convex function $\f:\R^n\to\R\cup\{+\infty\}$ and every $m\in\N^*$, we define the function $\f_m$ by 
\[
\f_m(x)=\left(\frac{m}{2}|\cdot|^2\right)\square\left(\f+\frac{1}{2m}|\cdot|^2\right)(x)=
\inf_{z}\left(\f(z)+\frac{|z|^2}{2m}+\frac{m}{2}|x-z|^2\right).
\]
Thus 
$$
\L\f_m(x)=\frac{|x|^2}{2m}+\L\f\square\left(\frac{m}{2}|\cdot|^2\right)(x)= \frac{|x|^2}{2m}+
\inf_z\left(\L\f(z)+\frac{m}{2}|z-x|^2\right).
$$
We shall need the following approximation lemma. 

\begin{lemma}\label{regul}
Let $n,m\ge1$ and $\f:\R^n\to\R\cup\{+\infty\}$ be even convex such that $0<\int_{\R^{n}}e^{-\f(x)}dx<+\infty$. Then, 
\begin{enumerate}
\item $\dom(\f_m)=\dom(\L\f_m)=\R^n$, $\f_m$ and $\L\f_m$ are differentiable and strictly convex on $\R^n$, $\nabla\f_m$ and $\nabla\L\f_m$ are bijective on $\R^n$. Moreover, $\nabla\f_m$ is $m$-Lipschitz continuous. 
\item When $m\to+\infty$ one has $\f_{m}(x)\to\f(x)$ and $\L\f_m(x)\to\L\f(x)$ a.e.
\item When $m\to +\infty$, for every measurable set $A$  one has $\int_Ae^{-\f_m}\to \int_Ae^{-\f}$, $\int_A\f_m e^{-\f_m}\to \int_A\f e^{-\f}$ and moreover for every $t>0$
\begin{eqnarray}\label{convergence}
P(t\f_m)\to P(t\f)\quad\hbox{and}\quad \frac{d}{dt}(P(t\f_m))\to\frac{d}{dt}(P(t\f)).
\end{eqnarray}
\end{enumerate}
\end{lemma}
\begin{proof}
1. For every $x\in\R^n$ one has $\f_m(x)\le\f(0)+\frac{m}{2}|x|^2<+\infty$, hence $\dom(\f_m)=\R^n$. In the same way $\L\f_m(x)\le\frac{|x|^2}{2m}+\L\f(0)+\frac{m}{2}|x|^2<+\infty$, hence $\dom(\L\f_m)=\R^n$. Moreover, it is clear that $\L\f_m$ is strictly convex and, using \cite{R} Theorem 26.3, it is not difficult to see that $\L\f_m$ is differentiable. It follows that the same holds for $\f_m$. The fact that $\nabla\f_m$  and $\nabla\L\f_m$ are bijective on $\R^n$ deduces from Corollary 26.3.1 in \cite{R}. Since $\L\f-\frac{1}{2m}|\cdot|^2$ is convex for some $m>0$ it follows from \cite{GR} that $\nabla\f$ is $m$-Lipschitz continuous.\\ 
2. These convergences are classical. Let us prove for example the first one. On one hand, one has $\f_m(x_0)\le\frac{|x_0|^2}{2m}+\f(x_0)$. On the other hand, if $x_0\in\dom(\f)$ then there exists a hyperplane touching the epigraph of $\f$ at $(x_0,\f(x_{0}))$, thus there exists $y\in\R^n$ such that $\f(x)\ge\f(x_0)+\langle x-x_0,y\rangle$, for all $x\in\R^n$. This implies that
\[
\f_m(x_0)\ge\inf_x \left(\f(x)+\frac{m}{2}|x-x_0|^2\right)\ge\inf_x \left(\f(x_0)+\langle x-x_0,y\rangle+\frac{m}{2}|x-x_0|^2\right)=\f(x_0)-\frac{|y|^2}{2m}.
\]
Letting $m\to+\infty$ gives the convergence. If $x_0\notin\overline{\dom(\f)}$ then, using that $\f\ge\min\f+I_{\dom(\f)}$, we deduce that
\[
\f_m(x_0)\ge\min\f+\inf_{x\in\dom(\f)}\frac{m}{2}|x-x_0|^2
=\min\f+\frac{m}{2}d(x_0,\dom(\f))^2.
\]
Therefore, $\f_{m}(x_0)\to\f(x_0)$ when $m\to+\infty$ for every $x_0\notin\partial(\dom(\f))$, that is a.e.
\\
3. First, notice that we may assume that $\f(0)=0$. From Lemma~\ref{lem:integ}, it follows that there exists $a,b>0$ such that, for every $x\in\R^n$, one has $a|x|-1\le\f(x)\le I_{bB_2^n}(x)+1$. Hence, we get
\begin{eqnarray}\label{eq:lower-bound-rho}
\f_m(x)\ge\inf_z \left(\f(z)+\frac{1}{2}|x-z|^2\right)\ge \inf_z \left(a|z|-1+\frac{1}{2}(|x|-|z|)^2\right)= a|x|-\frac{a^2}{2}-1.
\end{eqnarray}
Thus, $e^{-\f_{m}(x)} \leq e^{\frac{a^2}{2}+1-a|x|}$ for all $m$, then, from the dominated convergence theorem, one deduces that $\int_A e^{-\f_m} \to\int_A e^{-\f}$ when $m\to +\infty$. In the same way, one has $b|y|-1\le\L\f(y)\le I_{aB_2^n}(y)+1$ thus 
\[
\L\f_m(y)\ge\inf_z \left(\L\f(z)+\frac{1}{2}|y-z|^2\right)\ge b|y|-\frac{b^2}{2}-1.
\]
Hence, from the dominated convergence theorem one deduces that $\int_{A} e^{-\L\f_m}\to\int_{A} e^{-\L\f}$ when $m\to +\infty$. We conclude that $P(\f_m)\to P(\f)$ when $m\to +\infty$. 
Similarly, we prove that for every $t>0$ one has $P(t\f_{m})\to P(t\f)$ when $m\to +\infty$. 
Using that $ue^{-u}\leq\frac{2}{e}e^{-\frac{u}{2}}$ for every $u\in\R$ we get that $\f_{m}e^{-\f_{m}}\le \frac{2}{e}e^{-\frac{\f_{m}}{2}}$
and we conclude again by the dominated convergence theorem. The same method gives the result for $t\f_m$ and $\L(t\f_m)$.
\end{proof}

\section{General results on the functional volume product}

Let $T:\R^n\to\R^n$ be an invertible linear map. Then, putting $z=Tx$ we get, for every $y\in\R^n$,
\[
\L(\f\circ T)(y)=\sup_x(\langle x, y\rangle -\f(Tx))=\sup_z(\langle T^{-1}z, y\rangle -\f(z))=(\L\f)\left((T^{-1})^*(y)\right).
\]
Therefore, $\L(\f\circ T)=(\L\f)\circ (T^{-1})^*$. Hence, changing variables, we get $P(\f\circ T)=P(\f)$. The functional $P$ admits another invariance: for any $c\in\R$, one has $\L(\f+c)(y)=\L\f(y)-c$. Thus, 
\[
P(\f+c)=\int e^{-(\f+c)}\int e^{-\L\f+c}=P(\f).
\]
Hence, one may assume in the following that $\f(0)=0$. Since we are dealing with even functions, one has also $\f(0)=\min\f$ and $\L\f(0)=-\inf \f=-\f(0)$. Thus, if $\f(0)=0$ then $\L\f(0)=0$.
On the opposite, when $\f$ is replaced by $t\f$, for $t>0$, the functional $P$ is not invariant. We shall take advantage of this. For every $y\in\R^n$ and $t>0$, one has
\[
\L(t\f)(y)=\sup_x(\langle x, y\rangle -t\f(x))=t \sup_x(\langle x, \frac{y}{t}\rangle -\f(x))=t\L\f\left(\frac{y}{t}\right).
\]
Hence, changing variables, we get
\[
P(t\f)=\int_{\R^{n}}e^{-t\f(x)}dx\int_{\R^{n}}e^{-\L(t\f)(y)}dy=t^n\int_{\R^{n}}e^{-t\f(x)}dx\int_{\R^{n}}e^{-t\L\f(y)}dy.
\]
In the following, we denote by $\mu_\f$ the measure on $\R^n$ with density $e^{-\f}$ with respect to the Lebesgue measure and, for an oriented hypersurface $S$ of $\R^n$ whose normal $n_S$ is defined a.e., we denote
\[
V_{S}(\f)=\int_{S}n_{S}(y)e^{-\f(y)}dy\quad\hbox{and}\quad Q_S(\f)=\int_{S}\langle y,n_{S}(y)\rangle e^{-\f(y)}dy.
\]
Notice that, if $S$ is the boundary of a cone with apex at the origin, then $Q_S(\f)=0$.
The following proposition generalizes ideas from the proof of Theorem 10 in \cite{FM08a} and Proposition~1 in \cite{FHMRZ}. Let us recall that a subset $A$ of $\R^n$ is BV if its indicator function ${\bf 1}_A$ has bounded variation.

\begin{proposition}\label{prop:ineq-phi}
Let $n\ge1$ and $\f:\R^n\to\R\cup\{+\infty\}$ be differentiable and convex such that $0<\int_{\R^{n}}e^{-\f(x)}dx<+\infty$. Let $A$ be a bounded and BV subset of $\R^n$ such that $\mu_\f(A)>0$. Then, for any $x\in \R^n$ one has 
\begin{equation}\label{eq:ineq-phi}
\langle x,-\frac{V_{\partial A}(\f)}{\mu_\f(A)}\rangle-\f(x)\le n-\int_A\f(y)\frac{d\mu_\f(y)}{\mu_\f(A)}-\frac{Q_{\partial A}(\f)}{\mu_{\f}(A)},
\end{equation}
{\em i.e.} $\L\f\left(-\frac{V_{\partial A}(\f)}{\mu_\f(A)}\right)\le n-\int_A\f(y)\frac{d\mu_\f(y)}{\mu_\f(A)}-\frac{Q_{\partial A}(\f)}{\mu_{\f}(A)}$. Moreover, if for some $x_0\in\R^n$ there is equality in (\ref{eq:ineq-phi}) then $x_0\in\dom(\f)$ and $\f$ is affine on $[x_0,z]$ for every $z\in A\cap\dom(\f)$.
\end{proposition}

\begin{proof}
By convexity, one has, for all $y\in\dom(\f)$ and for all $x\in\R^n$,
\begin{equation*}
 \langle x,\nabla\f(y)\rangle-\f(x)\leq\langle y,\nabla\f(y)\rangle-\f(y). 
\end{equation*} 
We multiply by $e^{-\varphi(y)}$,  integrate in $y$ on $A$ and divide by $\mu_\f(A)$ to get 
\begin{equation}\label{eq:ineq-nabla}
\langle x,\int_A\nabla\f(y)\frac{d\mu_\f(y)}{\mu_\f(A)}\rangle-\f(x)\le\int_A\left(\langle y, \nabla\f(y)\rangle-\f(y)\right)\frac{d\mu_\f(y)}{\mu_\f(A)}.
\end{equation}
Recall the following consequence of the Stokes formula, known as Green's identities: for any sufficiently smooth $f,g:A\to\R$, one has 
\[
\int_A (f\Delta g+\langle \nabla f,\nabla g\rangle)=\int_{\partial A}f\langle \nabla g,n_A\rangle,
\]
where the integrals are taken with respect to the Hausdorff measure. Applying this formula to $f(y)=e^{-\varphi(y)}$ and $g(y)=\langle x,y\rangle$, where $x$ is a fixed vector, gives 
\begin{equation}\label{eq:green-nabla}
\int_A\nabla\f d\mu_\f=-\int_{\partial A}n_A(y)e^{-\f(y)}dy=-V_{\partial A}(\f).
\end{equation}
Applying it to $f(y)=e^{-\f(y)}$ and $g(y)=\frac{|y|^2}{2}$ gives
\begin{equation}\label{eq:green-nabla-Q}
\int_A\langle y, \nabla\f(y)\rangle d\mu_\f(y)=n\mu_\f(A)-Q_{\partial A}(\f).
\end{equation}
Substituting these values in the inequality (\ref{eq:ineq-nabla}), we conclude that inequality (\ref{eq:ineq-phi}) is verified. 
Moreover, if for some $x_0\in\R^n$ there is equality in (\ref{eq:ineq-phi}) then, for almost all $y\in\dom(\f)$, one has
$\varphi(y)+\langle x_0-y,\nabla\varphi(y)\rangle = \varphi(x_0)$. We conclude using Lemma 3 in \cite{FGMR}.
\end{proof}
 
\begin{remark}\label{rem1}
(i) Proposition \ref{prop:ineq-phi} and equalities (\ref{eq:green-nabla}) and (\ref{eq:green-nabla-Q}) hold if $\f$ is convex (without any differentiability assumption) and $A=\R_\e^n$, for some $\e\in\{-1;1\}^n$.\\
(ii) The inequality of Proposition \ref{prop:ineq-phi} holds if $\f$ is replaced by $\L\f$, $\L\f$ is differentiable, $\nabla\f$ is Lipschitz continuous and $A=\nabla\f(\R_\e^n)$, for some $\e\in\{-1;1\}^n$.
\end{remark}

\begin{proof}
(i) Let $R>0$ and $m\in\N$, we apply Proposition \ref{prop:ineq-phi} to $\f$ replaced by the function $\f_m$ introduced before Lemma~\ref{regul} and to $A_R=\R_\e^n\cap RB_2^n$.
\[
\langle x,-\frac{V_{\partial A_R}(\f_m)}{\mu_{\f_m}(A_R)}\rangle-\f_m(x)\le
n-\int_{A_R}\f_m(y)\frac{d\mu_{\f_m}(y)}{\mu_{\f_m}(A_R)}-\frac{Q_{\partial A_R}(\f_m)}{\mu_{\f_m}(A_R)},
\]
Then using Lemma~\ref{regul} and inequality (\ref{eq:lower-bound-rho}), when $m$ tends to $+\infty$, one has 
\[
\mu_{\f_m}(A_R)\to \mu_\f(A_R),\ \int_{A_R}\f_m d\mu_{\f_m}\to\int_{A_R}\f d\mu_{\f},\ V_{\partial A_R}(\f_m)\to V_{\partial A_R}(\f)\ \hbox{and}\ Q_{\partial A_R}(\f_m)\to Q_{\partial A_R}(\f).
\]
Now we let $R$ tend to $+\infty$. Clearly $\mu_\f(A_R)$ and $\int_{A_R}\f d\mu_{\f}$ converge respectively to $\mu_\f(A)$ and $\int_{A}\f d\mu_{\f}$. From Lemma \ref{lem:integ}, there exists $a>0$ such that $\f(x)\ge a|x|-1$, for every $x\in\R^n$. Thus, 
\[
\int_{R S^{n-1}\cap\R_\e^n}e^{-\f(x)}dx\le e^{-aR+1}R^{n-1}\frac{\omega_n}{2^n}\to 0, \quad\hbox{when}\ R\to+\infty,
\]
where $S^{n-1}$ is the Euclidean sphere of radius $1$ and $\omega_n$ is its $(n-1)$-dimensional measure. Hence, we easily see that $V_{\partial A_R}(\f)$ and $Q_{\partial A_R}(\f)$ converge respectively to $V_{\partial A}(\f)$ and $Q_{\partial A}(\f)$. In the same way, applying equalities (\ref{eq:green-nabla}) and (\ref{eq:green-nabla-Q}) to $\f_m$ and $A_R$, and letting $m$ and $R$ tend to infinity, we get that these equalities also hold for $\f$ and $A$. Moreover, if there is equality in (\ref{eq:ineq-phi}) then, using (\ref{eq:green-nabla}) and (\ref{eq:green-nabla-Q}), we get that there is equality in (\ref{eq:ineq-nabla}) and we conclude using Lemma 3 in \cite{FGMR}.\\
(ii) Let $R>0$ and $A_{R}=\nabla\f(\R^{n}_{\e}\cap RB^{n}_{2})$. Since $\nabla\f$ is Lipschitz continuous then $A_{R}$ is bounded and BV (see the remark p.119 in section 7.5 in \cite{Pf}). Thus we may apply Proposition \ref{prop:ineq-phi} to $A_R$ and $\f$ replaced by $\L\f$.  
In the same way as before, one sees that the terms appearing in the inequality converge when $R$ tends to $+\infty$.
\end{proof}

 \begin{corollary} \label{coro} 
 Let $n\ge1$ and $\f:\R^n\to\R\cup\{+\infty\}$ be differentiable and convex such that $0<\int_{\R^{n}}e^{-\f(x)}dx<+\infty$. 
 Let $A, B$ be bounded and BV subsets of $\R^n$ such that $\mu_\f(A),\mu_{\L\f}(B)>0$.
Then,
 \begin{equation}\label{ineq-phi-Lphi}
\langle \frac{V_{\partial A}(\f)}{\mu_\f(A)}, \frac{V_{\partial B}(\L\f)}{\mu_{\L\f}(B)}\rangle\le 2n-\int_A\f(y)\frac{d\mu_\f(y)}{\mu_\f(A)}-\int_B\L\f(y)\frac{d\mu_{\L\f}(y)}{\mu_{\L\f}(B)}-\frac{Q_{\partial A}(\f)}{\mu_\f(A)}-\frac{Q_{\partial B}(\L\f)}{\mu_{\L\f}(B)}.
\end{equation}
Moreover, if there is equality in (\ref{ineq-phi-Lphi}) then $\f$ is affine on $[-\frac{V_{\partial B}(\L\f)}{\mu_{\L\f}(B)},a]$ for every $a\in A\cap\dom(\f)$ and $\L\f$ is affine on $[-\frac{V_{\partial A}(\f)}{\mu_\f(A)},b]$ for every $b\in B\cap\dom(\L\f)$.
\end{corollary}

\begin{proof}
Since we are working with integrals, we may assume that $\f$ is lower semi-continuous. We apply the inequality $\langle x,y\rangle\le\f(x)+\L\f(y)$ to $x=-\frac{V_{\partial B}(\L\f)}{\mu_{\L\f}(B)}$ and $y=-\frac{V_{\partial A}(\f)}{\mu_{\f}(A)}$, apply Proposition \ref{prop:ineq-phi} twice to $\f$ and $\L\f$ and use that $\mathcal{L}(\mathcal{L}\varphi)=\varphi$ to deduce that
\begin{eqnarray*}
\langle \frac{-V_{\partial A}(\f)}{\mu_\f(A)}, \frac{-V_{\partial B}(\L\f)}{\mu_{\L\f}(B)}\rangle&\le&\L\f\left(-\frac{V_{\partial A}(\f)}{\mu_\f(A)}\right)+\f\left(-\frac{V_{\partial B}(\L\f)}{\mu_{\L\f}(B)}\right)\\
&\le&2n-\int_A\f(y)\frac{d\mu_\f(y)}{\mu_\f(A)}-\int_B\L\f(y)\frac{d\mu_{\L\f}(y)}{\mu_{\L\f}(B)}-\frac{Q_{\partial A}(\f)}{\mu_\f(A)}-\frac{Q_{\partial B}(\L\f)}{\mu_{\L\f}(B)}.
\end{eqnarray*}
Moreover, if there is equality in (\ref{ineq-phi-Lphi}) then there is equality in (\ref{eq:ineq-phi}) for $x=\frac{-V_{\partial B}(\L\f)}{\mu_{\L\f}(B)}$. Hence, from the equality case of Proposition \ref{prop:ineq-phi}, we deduce that $\f$ is affine on $[-\frac{V_{\partial B}(\L\f)}{\mu_{\L\f(B)}},a]$, for every $a\in A\cap\dom(\f)$. The same argument gives that $\L\f$ is affine on $[-\frac{V_{\partial A}(\f)}{\mu_\f(A)},b]$, for every $b\in B\cap\dom(\L\f)$.
\end{proof}
Notice that, changing variables, for $t>0$, one has
\[
\mu_{t\f}(tA)=\int_{tA}e^{-t\f(x)}dx=t^n\int_Ae^{-t\f(tz)}dz.
\]
Using again that $\L(t\f)(y)=t\L\f(\frac{y}{t})$, we have also
\[
\mu_{\L(t\f)}(tB)=\int_{tB}e^{-\L(t\f)(y)}dy=t^n\int_Be^{-t\L\f(z)}dz.
\]
We define 
\[
F_{A,B}(t)=\mu_{t\f}(tA)\mu_{\L(t\f)}(tB). 
\]
In the next lemma, we compute the derivatives of $\mu_{t\f}(tA), \mu_{\L(t\f)}(tB)$ and $F_{A,B}(t)$.

\begin{lemma}\label{cor:deriv}
 Let $n\ge1$ and $\f:\R^n\to\R\cup\{+\infty\}$ be differentiable and convex such that $0<\int_{\R^{n}}e^{-\f(x)}dx<+\infty$. 
 Let $A,B$ be bounded and BV subsets of $\R^n$ such that $\mu_\f(A),\mu_{\L\f}(B)>0$.
 Then,
 \begin{enumerate}
 \item ${\displaystyle (\mu_{t\f}(tA))'=-\int_{tA}\f d\mu_{t\f}+\frac{1}{t}Q_{t\partial A}(t\f)}$.
\item ${\displaystyle (\mu_{\L(t\f)}(tB))'=\frac{n}{t}\mu_{\L(t\f)}(tB)-\frac{1}{t}\int_{tB}\L(t\f) e^{-\L(t\f)}}$.
\item ${\displaystyle \frac{F_{A,B}'(t)}{F_{A,B}(t)}=\frac{n}{t}-\int_{tA}\f \frac{d\mu_{t\f}}{\mu_{t\f}(tA)}-\frac{1}{t}\int_{tB}\L(t\f) \frac{d\mu_{\L(t\f)}}{\mu_{\L(t\f)}(tB)}+\frac{1}{t}\frac{Q_{t\partial A}(t\f)}{\mu_{t\f}(tA)}}$.
\item In particular, if A is a cone with apex at the origin then
\[
F_{A,B}'(1)=nF_{A,B}(1)-\int_{A}\f e^{-\f}\int_{B}e^{-\L\f}-\int_{A}e^{-\f}\int_{B}\L\f e^{-\L\f}.
\]
\end{enumerate}
\end{lemma}

\begin{proof} 
 1. We compute the derivative of $\mu_{t\f}(tA)$, change variables and apply Green's identity (\ref{eq:green-nabla-Q}) to $t\f$ and $tA$, this gives
\begin{eqnarray*}
(\mu_{t\f}(tA))'&=&nt^{n-1}\int_Ae^{-t\f(tz)}dz-t^n\int_A(\f(tz)+t\langle\nabla\f(tz),z\rangle) e^{-t\f(tz)}dz\\
&=&\frac{n}{t}\mu_{t\f}(tA)-\int_{tA}\f d\mu_{t\f}-\int_{tA}\langle \nabla\f(x),x\rangle d\mu_{t\f}(x)\\
&=&-\int_{tA}\f d\mu_{t\f}+\frac{1}{t}Q_{t\partial A}(t\f).
\end{eqnarray*}
2., 3. and 4. The computation of the derivatives of $\mu_{\L(t\f)}(tB)$ and $F_{A,B}(t)$ are direct.
\end{proof}

 \begin{corollary} \label{coro-deriv} 
 Let $n\ge1$ and $\f:\R^n\to\R\cup\{+\infty\}$ be differentiable and convex such that $0<\int_{\R^{n}}e^{-\f(x)}dx<+\infty$. 
 Let $A, B$ be bounded and BV subsets of $\R^n$ such that $\mu_\f(A),\mu_{\L\f}(B)>0$.
 Then,
 \begin{equation}\label{ineq-deriv-F}
\frac{d}{dt}\left(t^nF_{A,B}\right)\ge  t^{n-1}\left(\langle V_{t\partial A}(t\f), V_{t\partial B}(\L(t\f))\rangle+2Q_{t\partial A}(t\f)\mu_{\L(t\f)}(tB)+Q_{t\partial B}(\L(t\f))\mu_{t\f}(tA)\right).
\end{equation}
Moreover, if there is equality in (\ref{ineq-deriv-F}) then $\f$ is affine on $[-\frac{V_{\partial B}(t\L\f)}{\mu_{t\L\f}(B)},a]$, for every $a\in (tA)\cap\dom(\f)$ and $\L\f$ is affine on $[-\frac{V_{\partial tA}(t\f)}{t\mu_{t\f}(tA)},b]$, for every $b\in B\cap\dom(\L\f)$.
\end{corollary}

\begin{proof}
Applying inequality (\ref{ineq-phi-Lphi}) of Corollary \ref{coro} to $tA$, $tB$ and $t\f$, and using 3) of Lemma~\ref{cor:deriv}  we get
\begin{eqnarray*}
\langle \frac{V_{t\partial A}(t\f)}{\mu_{t\f}(tA)}, \frac{V_{t\partial B}(\L(t\f))}{\mu_{\L(t\f)}(tB)}\rangle \le 
n+\frac{tF_{A,B}'(t)}{F_{A,B}(t)}-2\frac{Q_{t\partial A}(t\f)}{\mu_{t\f}(tA)}-\frac{Q_{t\partial B}(\L(t\f))}{\mu_{\L(t\f)}(tB)}.
\end{eqnarray*}
We multiply by $t^{n-1}F_{A,B}(t)$ and get inequality (\ref{ineq-deriv-F}).
\end{proof}

\begin{remark}\label{rk:asscoro}
Note that the inequalities of Corollary \ref{coro}, Lemma \ref{cor:deriv} and Corollary \ref{coro-deriv} hold in the two following cases:\\
(i) $\f$ and $\L\f$ convex without any differentiability assumption, $A=\R_\e^n$ and $B=T(\R_e^n)$, where $T$ is linear bijective and $\e\in\{-1;1\}^n$. Moreover, the equality cases also hold.\\
(ii) $\f$ and $\L\f$ are differentiable, $\nabla\f$ is Lipschitz continuous, $A=\R_\e^n$ and $B=\nabla\f(\R_e^n)$.
\end{remark}

\begin{lemma}\label{lem:V-R+}
Let $n\ge1$ and $\f:\R^n\to\R\cup\{+\infty\}$ be convex such that $0<\int_{\R^{n}}e^{-\f(x)}dx<+\infty$. Let $\varepsilon\in\{-1,1\}^n$. Then,\\
1) $\displaystyle{V_{\partial\R_\varepsilon^n}(\f)=-\sum_{i=1}^n \varepsilon_ie_i\int_{\R_\varepsilon^n\cap e_i^\bot}e^{-\f_i}}$.\\
2) If, moreover, $\f$ is differentiable, strictly convex, $\dom(\f)=\dom(\L\f)=\R^n$, $\nabla\f$ is Lipschitz continuous and the normal of $\nabla\f(\R_\varepsilon^n\cap e_i^\bot)$ is chosen exterior to  $\nabla\f(\R_\varepsilon^n)$ then 
\[ 
\langle V_{\nabla\f(\R_\varepsilon^n\cap e_i^\bot)}(\L\f),e_i\rangle=-\varepsilon_i\int_{\R_\varepsilon^n\cap e_i^\bot}e^{-\L(\f_i)}.
\]
\end{lemma}
\begin{proof}
1) It follows directly from the definition.\\
2) We assume that $\R_\varepsilon^n=\R_+^n$, the general case being the same. From the definition of $V$, one has
\[
\langle V_{\nabla\f(\R_+^n\cap e_i^\bot)}(\L\f),e_i\rangle=\int_{\nabla\f(\R_+^n\cap e_i^\bot)}\langle n_{\nabla\f(\R_+^n\cap e_{i}^{\perp})}(x),e_i\rangle e^{-\L\f(x)}dx.
\]
Now we use the notations and results obtained in Lemma \ref{Param}. Define 
$$
S_i=\nabla\f(\R_+^n\cap e_i^\bot)=\{y+t_i(y)e_i; y\in \R_+^n\cap e_i^\bot\}.$$ 
The surface $S_i$ is the graph of the function $t_i:\R_+^n\cap e_i^\bot\to\R$.
Each point $x\in S_i$ projects orthogonaly to $y=\nabla\f_i((\nabla\f)^{-1}(x))$.
Since $\nabla\f$ is Lipschitz continuous $S_i$ is rectifiable. Hence we get 
\[
\int_{S_i}\langle n_{\nabla\f(\R_+^n\cap e_i^\bot)}(y),e_i\rangle e^{-\L\f(y)}dy=-\int_{\R_+^n\cap e_i^\bot}e^{-\L\f(y+t_i(y)e_i)}dy.
\]
From Lemma \ref{Param}, we have $\L\f(y+t_i(y)e_i)=\L(\f_i)(y)$ for all $y\in e_i^\bot$ and thus we conclude. 
\end{proof}

\section{Proof in dimension $2$}
\begin{thm}\label{thm:dim2} Let $\f:\R^{2}\to\R\cup\{+\infty\}$ be even convex such that $0<\int_{\R^{2}}e^{-\f(x)}dx<+\infty$, then 
\[
P(\f)=\int_{\R^{2}}e^{-\f(x)}dx\int_{\R^{2}}e^{-\L\f(y)}dy\geq 4^{2}=16.
\]
Moreover, if $\f$ is lower semi-continuous, then $P(\f)=16$ if and only if there exists $a\in\R$ such that, either $\f=I_P+a$, or $\f=\|\cdot\|_{P}+a$, with $P$ being a parallelogram centered at the origin or there exists a basis $(u_1,u_2)$ of $\R^2$ and $b,c>0$ such that $\f(x_1u_1+x_2u_2)=c|x_1|+I_{[-b,b]}(x_2)+a$ for every $x_1,x_2\in\R$.
\end{thm}

\subsection{The inequality in dimension $2$}

Let us give an outline of the proof of the inequality of Theorem \ref{thm:dim2}.
We first show, in Lemma \ref{lem:equip}, that one can reduce to the case where $\f$ is strongly equipartitioned, in the sense that
$$
\f(0)=0,\quad \int_{\R_+}e^{-\f(te_{1})}dt=\int_{\R_+}e^{-\f(te_{2})}dt=1, \quad\int_{\R^{2}_{+}}e^{-\f}=\frac{1}{4}\int_{\R^{2}}e^{-\f}, \quad\int_{\R^{2}_{+}}\f e^{-\f}=\frac{1}{4}\int_{\R^{2}}\f e^{-\f}.
$$
Then, using Lemma \ref{regul}, we prove that we may assume that $\f$ has some nice regularity properties. From the strong equipartition of $\f$ and Lemma~\ref{cor:deriv}, we deduce that 
\[
P(\f)=8(F_1(1)+F_2(1))\quad \hbox{and}\quad \frac{d}{dt}\left(t^{2}P(t\f)\right)_{|t=1}=8\frac{d}{dt}(t^{2}(F_{1}(t)+F_{2}(t)))_{|t=1},
\]
where $F_1(t)=F_{\R_+^2,\nabla\f(\R_+^2)}(t)$ and $F_2(t)=F_{\R_+\times\R_{-},\nabla\f(\R_+\times\R_{-})}(t)$. 
Using Corollary \ref{coro-deriv} and Lemma \ref{lem:V-R+}, we show that 
\[
\frac{d}{dt}\left(t^{2}P(t\f)\right)_{|t=1}\ge 16\left(\int_0^{+\infty}e^{-\L(\f_1)(y)}dy+\int_0^{+\infty}e^{-\L(\f_2)(y)}dy\right).
\]
The inequality on dimension 1, established in \cite{FGMR}, shows that $\frac{d}{dt}\left(t^{2}P(t\f)\right)_{|t=1}\ge 32$. Then, we conclude by integrating.

\begin{lemma}\label{lem:equip}
Let $\f:\R^{2}\to\R\cup\{+\infty\}$ be an even convex function such that $0<\int_{\R^{2}}e^{-\f(x)}dx<+\infty$. Then there exists a linear invertible map $T$ such that the function $\Tilde{\f}$ defined by $\Tilde{\f}(x)=\f\circ T(x)-\f(0)$ is strongly equipartitioned. 
\end{lemma}
\begin{proof}[\bf{Proof}]
For any $u\in S^1$, let $C(u)\subset S^1$ be the open half-circle delimited by $u$ and $-u$ containing the vectors $v$ which are after $u$ with respect to the counterclockwise orientation of $S^1$. For $v\in C(u)$, let  $C_{u,v}=\R_+u+\R_+v$ be the cone generated by $u$ and $v$ and define $f_u(v)=\mu_\f(C_{u,v})$. The map $f_u$ is continuous and increasing on $C(u)$, $f_u(u)=0$ and, since $\f$ is even,  $f_u(v)\longrightarrow\mu_{\f}(\R^{2})/2$ when $v\longrightarrow -u$. Thus, there exists a unique $v(u)\in C(u)$ such that $f_u(v(u))=\mu_\f(C_{u,v(u)})=\mu_\f(\R^2)/4$. Notice that $v:S^1\to S^1$ is continuous and, since $\f$ is even, one has $v(v(u))=-u$ for any $u\in S^1$.
For $u\in S^1$, let  $g(u)=\int_{C_{u,v(u)}}\f e^{-\f}-\frac{1}{4}\int_{\R^{2}}\f e^{-\f}$. Then, $g$ is continuous on $S^1$ and, since $\f$ is even, 
$$
g(u)+g(v(u))=\int_{C_{u,v(u)}}\f e^{-\f}+\int_{C_{v(u),-u}}\f e^{-\f}-\frac{1}{2}\int_{\R^{2}}\f e^{-\f}=0.
$$
Hence, $g(u)=-g(v(u))$. By the intermediate value theorem, there exists $u\in S^1$ such that $g(u)=0$, thus 
\[\int_{C_{u,v(u)}}\f e^{-\f}=\frac{1}{4}\int_{\R^{2}}\f e^{-\f}\quad\hbox{and}\quad \mu_\f\left(C_{u,v(u)}\right)=\frac{1}{4}\mu_\f\left(\R^{2}\right).
\]
Let $S$ be the linear map defined by $S(e_{1})=u$ and $S(e_{2})=v(u)$, then  $S(\R^{2}_{+})=C_{u,v(u)}$. Moreover, changing variables, for any Borel set $A$ in $\R^2$ we have $\mu_{\f\circ S}(A)=\mu_\f(S(A))/\det(S)$, thus,
$$
\mu_{\f\circ S}(\R^{2}_{+})=\frac{\mu_{\f}(S(\R^{2}_{+}))}{\det(S)}=\frac{\mu_{\f}(C_{u,v(u)})}{\det(S)}=\frac{\mu_{\f}(\R^{2})}{4\det(S)}=\frac{\mu_{\f\circ S}(\R^{2})}{4}.
$$
In the same way, one has 
\[
\int_{\R_+^2}(\f\circ S)e^{-\f\circ S}=\frac{1}{4}\int_{\R^2}(\f\circ S)e^{-\f\circ S}.
\]
Let $\alpha_{i}=\int_{0}^{+\infty}e^{-\f(re_{i})}dr$ and $\Delta$ be the linear map defined by $\Delta(e_{i})=\alpha_{i}e_{i}$ and $T=S\circ\Delta$, then a change of variables shows that $\Tilde{\f}=\f\circ T-\f(0)$ is strongly equipartitioned. 
\end{proof}

\begin{proof}[\bf{Proof of the inequality of Theorem \ref{thm:dim2}}]
Let $\f:\R^{2}\to\R\cup\{+\infty\}$ be an even convex function such that $0<\int_{\R^{2}}e^{-\f(x)}dx<+\infty$. From Lemma \ref{lem:equip}, there exists a linear invertible map $T$ such that the function $\Tilde{\f}$ defined by $\Tilde{\f}(x)=\f\circ T(x)-\f(0)$ is strongly equipartitioned. 
Since $P(\f)=P(\Tilde\f)$, we may assume that $\varphi$ is strongly equipartitioned. 
From Lemma \ref{regul}, we can assume that $\dom(\f)=\dom(\L\f)=\R^2$, $\f$ is differentiable and strictly convex on $\R^2$.  Then, up to sets of Lebesgue measure zero, we have the partition $\R^2=\cup_{\varepsilon\in\{-1,1\}^{2}}\nabla\f(\R_{\varepsilon}^{2})$. Using the equipartition , we get
$$
P(\f)=\sum_{\varepsilon\in\{-1,1\}^{2}}\int_{\R^{2}}e^{-\f}\int_{\nabla\f(\R_{\varepsilon}^{2})}e^{-\mathcal{L}\f}=4\sum_{\varepsilon\in\{-1,1\}^{2}}\int_{\R_{\varepsilon}^{2}}e^{-\f}\int_{\nabla\f(\R_{\varepsilon}^{2})}e^{-\mathcal{L}\f}.
$$ 
Using the fact that $\f$ is even, we get 
\[
P(\f)=8\int_{\R_+^{2}}e^{-\f}\int_{\nabla\f(\R_+^{2})}e^{-\L\f}+8\int_{\R_+\times\R_{-}}e^{-\f}\int_{\nabla\f(\R_+\times\R_{-})}e^{-\L\f}=8(F_{1}(1)+F_{2}(1))
\]
where $F_1(t)=F_{\R_+^2,\nabla\f(\R_+^2)}(t)$ and $F_2(t)=F_{\R_+\times\R_{-},\nabla\f(\R_+\times\R_{-})}(t)$. Using 4) of Lemma \ref{cor:deriv} and Remark~\ref{rk:asscoro}, we have 
$$
F_{1}'(1)=2F_{1}(1)-\int_{\R^{2}_{+}}\f e^{-\f}\int_{\nabla\f(\R^{2}_{+})}e^{-\L\f}-\int_{R^{2}_{+}}e^{-\f}\int_{\nabla\f(\R^{2}_{+})}\L\f e^{-\L\f}.
$$
$$
F_{2}'(1)=2F_{2}(1)-\int_{\R_{+}\times\R_{-}}\f e^{-\f}\int_{\nabla\f(\R_{+}\times\R_{-})}e^{-\L\f}-\int_{\R_{+}\times\R_{-}}e^{-\f}\int_{\nabla\f(\R_{+}\times\R_{-})}\L\f e^{-\L\f}.
$$
Then, we get
\begin{equation*}\label{equation1}
\begin{aligned}
\frac{d}{dt}(t^{2}(F_{1}(t)+F_{2}(t)))_{|t=1}&=4(F_{1}(1)+F_{2}(1))-\left(\int_{\R^{2}_{+}}\f e^{-\f}\int_{\nabla\f(\R^{2}_{+})}e^{-\L\f}+\int_{\R^{2}_{+}}e^{-\f}\int_{\nabla\f(\R^{2}_{+})}\L\f e^{-\L\f}\right.\\
&+\left.\int_{\R_{+}\times\R_{-}}\f e^{-\f}\int_{\nabla\f(\R_{+}\times\R_{-})}e^{-\L\f}+\int_{\R_{+}\times\R_{-}}e^{-\f}\int_{\nabla\f(\R_{+}\times\R_{-})}\L\f e^{-\L\f}\right).
\end{aligned}
\end{equation*}
Thus, using the fact that $\f$ is equipartitioned, we get that
\begin{equation*}
\begin{aligned}
\frac{d}{dt}\left( t^{2}(F_{1}(t)+F_{2}(t)\right))_{|t=1}&=4(F_{1}(1)+F_{2}(1))-\int_{\R^{2}_{+}}\f e^{-\f}\times\frac{1}{2}\int_{\R^{2}}e^{-\L\f}-\int_{\R^{2}_{+}}e^{-\f}\times\frac{1}{2}\int_{\R^{2}}\L\f e^{-\L\f}\\
&=\frac{1}{2}P(\f)-\frac{1}{8}\left(\int_{\R^{2}}\f e^{-\f}\int_{\R^{2}}e^{-\L\f}+\int_{\R^{2}}e^{-\f}\int_{\R^{2}}\L\f e^{-\L\f} \right).
\end{aligned}
\end{equation*}
On the other hand, applying Lemma \ref{cor:deriv} and Remark~\ref{rk:asscoro} 
one has
\begin{eqnarray*}
\frac{d}{dt}\left(t^{2}P(t\f)\right)_{|t=1}
&=&4P(\f)-\left(\int_{\R^{2}}\f e^{-\f}\int_{\R^{2}}e^{-\L\f}+\int_{\R^{2}}e^{-\f}\int_{\R^{2}}\L\f e^{-\L\f}\right)\\&=&8\frac{d}{dt}(t^{2}(F_{1}(t)+F_{2}(t)))_{|t=1}.
\end{eqnarray*}
We apply Corollary \ref{coro-deriv} and Remark~\ref{rk:asscoro} for $A=\R_+^{2}$ and $B=\nabla \f(\R_+^{2})$ and use the equipartition to get 
\begin{eqnarray}\label{ineq:++}
\frac{d}{dt}\left(t^2F_1(t)\right)_{|t=1}\ge  \langle V_{\partial \R_+^2}(\f), V_{\partial \nabla \f(\R_+^{2})}(\L\f)\rangle+\frac{\mu_{\f}(\R^{2})}{4}Q_{\partial \nabla \f(\R_+^{2})}(\L\f).
\end{eqnarray} 
From Lemma \ref{lem:V-R+}, one has 
\begin{eqnarray}\label{eq:V-dim2}
V_{\partial\R^{2}_{+}}(\f)=-e_1\int_0^{+\infty}e^{-\f_1}-e_2\int_0^{+\infty}e^{-\f_2}=-(e_{1}+e_{2}).
\end{eqnarray}
Similarly, one has $V_{\partial \nabla \f(\R_+^{2})}(\L\f)=-W_1-W_2$, where 
\[
W_1=-\int_{\nabla\f(\{0\}\times\R_+)}n_{\nabla\f(\R_+^2)}e^{-\L\f}
\quad\hbox{and}\quad
W_2=-\int_{\nabla\f(\R_+\times\{0\})}n_{\nabla\f(\R_+^2)}e^{-\L\f}.
\]
Thus, the equation (\ref{ineq:++}) becomes 
\begin{eqnarray}\label{ineq:F1}
\frac{d}{dt}\left(t^2F_1(t)\right)_{|t=1}\ge \langle e_{1}+e_{2},W_1+W_2\rangle+\frac{\mu_{\f}(\R^{2})}{4}Q_{\partial \nabla \f(\R_+^{2})}(\L\f).
\end{eqnarray} 
Moreover, since $\f$ is even, one has 
\[
V_{\partial(\R_+\times\R_{-})}(\f)=-e_{1}+e_{2},\quad V_{\partial\nabla\f(\R_+\times\R_{-})}(\L\f)=-W_1+W_2
\]
and
\[
Q_{\partial \nabla \f(\R_+\times\R_{-})}(\L\f)=\int_{\partial \nabla \f(\R_+\times\R_{-})}\langle y,n_{\nabla \f(\R_+\times\R_{-})}(y)\rangle e^{-\L\f(y)}dy=-Q_{\partial \nabla \f(\R_+^{2})}(\L\f).
\]
Applying  Corollary \ref{coro-deriv} and Remark~\ref{rk:asscoro} for $A=\R_+\times\R_{-}$ and $B=\nabla \f(\R_+\times\R_{-})$ and using the equipartition, we get
\begin{eqnarray}\label{ineq:F2}
\frac{d}{dt}\left(t^2F_2(t)\right)_{|t=1}\ge \langle e_1-e_2,W_1-W_2\rangle-\frac{\mu_{\f}(\R^{2})}{4}Q_{\partial \nabla \f(\R_+^{2})}(\L\f).
\end{eqnarray} 
Adding (\ref{ineq:F1}) and (\ref{ineq:F2}), 
we obtain 
\[
\frac{d}{dt}\left(t^2P(t\f)\right)_{|t=1}= 8\frac{d}{dt}(t^{2}(F_{1}(t)+F_{2}(t)))_{|t=1}
\ge 16(\langle e_1,W_1\rangle +\langle e_2,W_2\rangle).
\]
Moreover, from Lemma \ref{lem:V-R+} for $i=1,2$, one has 
\begin{eqnarray*}
\langle e_i,W_i\rangle=\int_{\nabla\f(\R^{2}_{+}\cap e_{i}^{\perp})}\langle -e_i,n_{\nabla\f(\R_+^2)}(x)\rangle e^{-\L\f(x)}dx=\int_0^{+\infty}e^{-\L(\f_i)(y)}dy\ge 1,
\end{eqnarray*}
where the last inequality comes from the result in dimension 1 proved in \cite{FM08a, FM08b}. Thus, we get 
\begin{eqnarray}\label{ineq:final-phi-i}
\frac{d}{dt}\left(t^2P(t\f)\right)_{|t=1}\ge 16\left( \int_0^{+\infty}e^{-\L(\f_1)}+\int_0^{+\infty}e^{-\L(\f_2)}\right)\ge 32.
\end{eqnarray}
Applying this relation for $\f$ replaced by $s\f$ and using the fact that
\begin{eqnarray*}
\frac{d}{dt}\left(t^{2}P(t s\f)\right)_{|t=1}&=&\lim_{t\to 1}\frac{t^{2}P(ts\f)-P(s\f)}{t-1}=\lim_{u\to s}\frac{\frac{u^{2}}{s^{2}}P(u\f)-P(s\f)}{\frac{u}{s}-1}=\frac{1}{s}\frac{d}{du}\left(u^{2}P(u\f)\right)_{|u=s}
\end{eqnarray*} 
one gets, $\forall t>0$, $\frac{d}{dt}(t^{2}P(t\f))\geq 32t$.
Integrating this inequality, we conclude that for every $0<\e<t$
\begin{eqnarray}\label{ineq:final-int}
t^2P(t\f)= \e^2P(\e\f)+\int_\e^t\frac{d}{ds}\left(s^2P(s\f)\right)ds\ge 32\int_\e^t sds=16(t^2-\e^2).
\end{eqnarray}
Letting $\e$ tends to $0$, we conclude that $t^2P(t\f)\ge 16t^2$ and thus $P(\f)\ge16$.$\hfill\blacksquare$ \\

\subsection{The equality case in dimension $2$}
Now, we establish the equality case. 
We first prove, in Lemmas \ref{lem:equip-dim1} and \ref{lem:bounds}, that the set of even convex functions, which are equipartitioned, has some compactness property. Then, using this compactness, we show that, for any even convex function $\f$, there exists a sequence of smooth even convex strongly equipartitioned functions $\psi_k$ and an invertible linear transform $T$ such that $\psi_k\to \f\circ T$. Using this approximation, we show that the inequalities (\ref{ineq:final-phi-i}) hold for $\f$,  without any regularity assumption. We deduce that the restrictions $\f_1$ and $\f_2$ of $\f$ satisfy the equality case in dimension~1. From \cite{FGMR} it follows that either $\f_i=I_{[-1,1]}$, or $\f_i(x)=|x|$, for $i=1,2$. Then, we discuss three cases and we conclude.

We adapt the arguments of \cite{FHMRZ} to the functional case. This requires to develop new functional inequalities.


\begin{lemma}\label{lem:equip-dim1}
Let  $\f:\R_+\to\R\cup\{+\infty\}$ be a convex, non-decreasing function such that $\f(0)=0$ and $\int_0^{+\infty}e^{-\f(x)}dx=1$. Then, for every $x\in\R_+$
 \[
 x-1\le\f(x)\le I_{[0,1]}(x)+x\le I_{[0,1]}(x)+1.
 \]
 \end{lemma}

\begin{proof}
For $t\in[0,1]$, we define  $\psi(t)=\sup\{x\ge0; \f(x)\le t\}$. One has $\f(\psi(t))\le t$ for almost all $t\ge0$. From Jensen's inequality, we get
\[
\f(1)=\f\left(\int_0^{+\infty}e^{-\f(x)}dx\right)= \f\left(\int_0^{+\infty} \psi(t)e^{-t}dt\right)\le\int_0^{+\infty} \f(\psi(t))e^{-t}dt\le \int_0^{+\infty} te^{-t}dt=1.
\]
By convexity, we deduce that $\f(x)=\f((1-x)\cdot 0+x\cdot 1)\le x$ for $x\in[0,1]$ and thus $\f(x)\le I_{[0,1]}(x)+x\le I_{[0,1]}(x)+1$, for every $x\ge0$.
For the proof of the lower bound, the idea is exactly the same as in the proof of Proposition \ref{prop:ineq-phi}. 
By convexity, the function $\varphi$ is differentiable almost everywhere on $\dom(\f)$ and one has, for almost all $y\in\dom(\f)$ and for all $x\ge0$,
\begin{equation*}
\f(x)\ge \varphi(y)+ (x-y)\f'(y).
\end{equation*} 
We multiply by $e^{-\varphi(y)}$ and integrate in $y$ on $[0,+\infty)$ to get that, for every $x\ge0$,
\[
\f(x)\ge \int_0^{+\infty}\f(y) e^{-\f(y)}dy+\int_0^{+\infty}(x-y)\f'(y)e^{-\f(y)}dy.
\]
Using that $\f\ge0$ and  integrating by parts, we deduce that, for every $x\ge0$,
\[
\f(x)\ge\int_0^{+\infty}(x-y)\f'(y)e^{-\f(y)}dy=x-1.
\]
\end{proof}
Now, we prove the analogue lemma in dimension 2. Recall that, for $x=(x_1,\dots,x_n)\in\R^n$, one denotes $\|x\|_1=\sum_{i=1}^n|x_i|$ and $B_1^n=\{x\in\R^n;\|x\|_1\le1\}$.
\begin{lemma}\label{lem:bounds}
Let $\f:\R^2\to\R\cup\{+\infty\}$ be a convex, even and strongly equipartitioned function such that $\f(0)=0$. Then, for every $x\in\R^2$, one has
\[
\frac{\|x\|_1}{e+2}-2\le \f(x)\le I_{B_1^2}(x)+1 \quad\hbox{thus}\quad \frac{2}{e}\le \int_{\R^2}e^{-\f(x)}dx\le \left(2e(e+2)\right)^2.
\]
\end{lemma}

\begin{proof}
The bounds on the integral follows directly from the bounds on the function. Let us prove these bounds. 
From Lemma \ref{lem:equip-dim1} applied to $t\mapsto \f(te_i)$, we deduce that $\f(e_i)\le 1$, for $1\le i\le 2$. Since $\f$ is even and convex, we deduce that $\f(x)\le1$, for every $x\in B_1^2=\conv(\pm e_1, \pm e_2)$. This proves the upper bound.  

To prove the lower bound, define $c=(e+2)^{-1}<1$.
Let assume by contradiction that there exists $a=(a_1,a_2)\in\R^2$ such that 
\[
\f(a)<c\|a\|_1-2.
\]
By symmetry, we may assume that $a_1\ge a_2\ge0$. Moreover, if $a_2=0$, from Lemma~\ref{lem:equip-dim1} applied to $\varphi_2$, one has $\f(a)=\f_2(a_1)\ge |a_1|-1= \|a\|_1-1$, which is not possible since $c<1$. Thus, one has $a_2>0$. For every $x=(x_1,x_2)\in\R_+\times\R_-$,  by convexity, one has
\[
\f\left(\frac{a_2}{a_2-x_2}x+\frac{-x_2}{a_2-x_2}a\right)\le\frac{a_2}{a_2-x_2}\f(x)+\frac{-x_2}{a_2-x_2}\f(a).
\]
Since $a_2x-x_2a=(a_2x_1-x_2a_1)e_1\in \R_+e_1$, we may apply Lemma~\ref{lem:equip-dim1} and get
\[
\f\left(\frac{a_2}{a_2-x_2}x+\frac{-x_2}{a_2-x_2}a\right)\ge \frac{a_2}{a_2-x_2}x_1+\frac{-x_2}{a_2-x_2}a_1-1.
\]
Thus, we deduce that
\[
\f(x)\ge x_1+\frac{(-x_2)}{a_2}\left(-\f(a)-1+a_1\right)-1.
\]
Since $a_1\geq a_2>0$, using the upper bound on $\f(a)$, we get 
\[
-\f(a)-1+a_1\ge-c(a_1+a_2)+a_1+1\ge(1-c)a_1-ca_2\ge (1-2c)a_2.
\]
Therefore, for every $x\in\R_+\times\R_-$, we deduce that
\[
\f(x)\ge x_1+(1-2c)(-x_2)-1.
\]
Integrating on $\R_+\times\R_-$ and replacing $c$ by its value, we get 
\[
\int_{\R_+\times\R_-}e^{-\f(x)}dx\le \frac{e}{1-2c}=\frac{1}{c}.
\]
Now, we apply Proposition \ref{prop:ineq-phi} to the cone $A=\R_+^2$. Recall that, in this case, the term $Q_{\partial A}(\f)$ vanishes. Thus, we get
\[
\f(a)\geq -2+\langle a,-\frac{V_{\partial\R_+^2}(\f)}{\mu_{\f}(\R_+^2)}\rangle+\int_{\R_+^2}\f(y)\frac{d\mu_{\f}(y)}{\mu_{\f}(\R_+^2)}.
\]
Using that $\f\ge0$ and $V_{\partial\R_+^2}(\f)=-(e_1+e_2)$, that was established in equation (\ref{eq:V-dim2}), we get
\[
\f(a)\ge -2+\langle a, \frac{e_{1}+e_{2}}{\mu_\f(\R_+^2)}\rangle=-2+\frac{\|a\|_1}{\mu_\f(\R_+^2)}\ge-2+c\|a\|_1,
\]
which is in contradiction with the definition of $c$.
\end{proof}
\begin{remark}\label{rem:value:cn}
The same argument applies inductively and shows that there exists a constant $c_n>0$ such that, for any even equipartitioned convex function $\f$ on $\R^n$, such that  all its restrictions to coordinate subspaces of all dimensions are equipartitioned and $\f(0)=0$, one has
\[
c_n\|x\|_1-n\le \f(x)\le I_{B_1^n}(x)+1 \quad\hbox{thus}\quad \frac{2^n}{en!}\le \int_{\R^n}e^{-\f(x)}dx\le \left(\frac{2e}{c_n}\right)^n.
\]
\end{remark}

\begin{lemma}\label{lem:eq:conv}
Let $\f:\R^2\to\R\cup\{+\infty\}$ be an even lower semi-continuous convex function such that $0<\int e^{-\f}<+\infty$ and $\f(0)=0$. Then, there exists a sequence $(\psi_k)_k$ of differentiable strongly equipartitioned even strictly convex functions with $\dom(\psi_k)=\R^2$, such that $\nabla\psi_k$ is Lipschitz continuous and there exists an invertible linear map $T$ such that \\
(i) for every $x\in\R^2$, $(\psi_k(x))_k$ converges to $\f\circ T(x)$ and  $e^{-\psi_k(x)}\le Ce^{-d|x|}$, for some $C,d>0$\\
(ii) for every $x\in\R^2$, $(\L\psi_k(x))_k$ converges to $\L(\f\circ T)(x)$  and $e^{-\L\psi_k(x)}\le Ce^{-d|x|}$, for some $C,d>0$. 
\end{lemma}

\begin{proof}
We define the set $K_\f=\{x\in\R^2; \f(x)\le 1\}$. 
From Lemma \ref{lem:integ}, $K_\f$ is a symmetric convex body, 
and there exists $a,b>0$ such that $a|x|-1\le\f(x)\le I_{bB_2^2}(x)+1$. We recall the definition of $\f_m$ used in Lemma~\ref{regul}: for every $x\in\R^2$,
\[
\f_m(x)=\inf_z \left(\f(z)+\frac{|z|^2}{2m}+\frac{m}{2}|x-z|^2\right).
\]
Using the lower bound obtained in (\ref{eq:lower-bound-rho}), we have $\f_m(x)\ge a|x|-\frac{a^2}{2}-1$.
Hence, for every $m\in\N^{*}$, 
\[\{x;\f_m(x)\le1\}\subset RB_2^2,\]
where $R=\frac{a}{2}+\frac{2}{a}$.
There exists a sequence of invertible linear maps $T_m$ such that $\f_m\circ T_m$ is strongly equipartitioned. From Lemma~\ref{lem:bounds}, one has $\f_m(T_m(e_i))\le1$, for all $1\le i\le 2$, thus 
\[
T_m(B_1^2)\subset \{x;\f_m(x)\le1\}\subset RB_2^2.
\]
Hence, the sequence $(T_m)_m$ is bounded in the normed spaces of linear maps and thus there exists a subsequence $(T_{m_k})_k$ of linear maps that converges to some linear map $T$. Let us prove that $T$ is invertible. For every $m\in\N^{*}$, using Lemma~\ref{lem:bounds} and denoting $c=(e+2)^{-1}$, one has $c\|x\|_1-2\le\f_m(T_mx)$ for every $x$. Moreover, since $\f_m(x)\le\f(x)+\frac{|x|^2}{2}$ and $\f(x)\le I_{bB_2^2}(x)+1$, it follows that
\[
\f_m(x)\le I_{bB_2^2}(x)+1+\frac{|x|^2}{2}\le I_{bB_2^2}(x)+ 1+\frac{b^2}{2}.
\]
Thus, for any $x\in\R^2$,
\[
\frac{bc\|x\|_1}{|T_mx|}\le\f_m\left(\frac{bT_mx}{|T_mx|}\right)+2\le 3+\frac{b^2}{2}.
\]
This gives that, for every $x\in\R^2$,
\[
\left(\frac{3}{b}+\frac{b}{2}\right)|T_mx|\ge c\|x\|_1.
\]
Hence, $T$ satisfies the same bound and thus is invertible. Moreover, since $\f_m(T_mx)\ge c\|x\|_1-2$, we conclude that the sequence $\psi_k=\f_{m_k}\circ T_{m_k}$ is a sequence of strongly equipartitioned even differentiable strictly convex functions such that $(\psi_k(x))_k$ converges to $\f\circ T(x)$ and $e^{-\psi_k(x)}\le e^{2-c\|x\|_1}$. Thus, $\f\circ T$ is strongly equipartitioned. Moreover, the Lipschitz continuity of $\nabla\psi_k$ is clearly inherited from the one of $\nabla\f_m$ and from Lemma~\ref{lem:bounds}, one has $\f_m(T_mx)\le I_{B_1^2}(x)+1$, hence
\[
\L(\f_m\circ T_m)(x)\geq\L\left(I_{B_1^2}+1\right)(x)= \|x\|_\infty-1.
\]
Therefore, $\L\psi_k=\L(\f_{m_k}\circ T_{m_k})$ satisfies the same bound. From Lemma~\ref{regul}, one has $\L\f_m(x)\to\L\f(x)$, for every $x\in\R^2$, when $m\to+\infty$. Since $T_{m_k}$ converges to $T$, we conclude that, for every $x\in\R^{2}$, $\L\psi_k(x)$ converges to $\L(\f\circ T)(x)$.
\end{proof}

\noindent
{\bf{Proof of the equality case.}} Let $\f:\R^{2}\to\R\cup\{+\infty\}$ be an even lower semi-continuous convex function such that $0<\int_{\R^{2}}e^{-\f(x)}dx<+\infty$ and $P(\f)=16$. By Lemma \ref{lem:eq:conv}, there exists a sequence $(\psi_k)_k$ of differentiable strongly equipartitioned even strictly convex functions with $\dom(\psi_k)=\R^2$ and a bijective linear map $T$ such that $\f\circ T$ is strongly equipartitioned and $i)$ and $ii)$ of Lemma \ref{lem:eq:conv} are satisfied. Since $P(\f\circ T)=P(\f)$ and our equality case is invariant by invertible linear maps, we replace $\f\circ T$ by $\f$ in the rest of the proof.  
We have thus established that $\f$ is the limit of a sequence $(\psi_k)_k$ of differentiable  and strongly equipartitioned strictly convex functions. Thus, the inequalities (\ref{ineq:final-phi-i}) and (\ref{ineq:final-int}) are valid for the functions $\psi_k$. From $i)$ and $ii)$ of Lemma~\ref{lem:eq:conv}, we may apply the dominated convergence theorem and, taking the limit, using Lemma~\ref{regul}, we deduce that these inequalities are also valid for $\f$.  Then, applying the same reasoning  as in inequality (\ref{ineq:final-int}) for $t=1$, and using that $P(\f)=16$ and  $P(\e\f)\ge16$, we get
\[
16=P(\f)=\e^2P(\e\f)+\int_\e^1\frac{d}{ds}\left(s^2P(s\f)\right)ds\ge 16\e^2 +32\int_\e^1 sds=16\e^2+16(1-\e^2)=16.
\]
Thus, there is equality in the intermediate inequalities. Hence, for every $0<\e\le1$, one has $P(\e\f)=16$. Thus, for all $0<t\le1$, we have $\frac{d}{dt}\left(t^2P(t\f)\right)=\frac{d}{dt}\left(16t^2\right)=32t$. Hence, there is equality in (\ref{ineq:final-phi-i}). This implies that 
\[
\int_0^{+\infty}e^{-\L(\f_1)}=1\quad\hbox{and}\quad\int_0^{+\infty}e^{-\L(\f_2)}=1.
\] 
From the equality case in dimension $1$, due to \cite{FGMR}, and since $\f_{i}$ is equipartitioned, we deduce that, for $i=1,2$, either $\f_i(x)=I_{[-1,1]}(x)$, or $\f_i(x)=|x|$, for every $x\in\R$. Following the proof in \cite{FGMR}, we distinguish three cases. \\

{\bf A.} If $\f_2(x)=\f(x,0)=I_{[-1,1]}(x)$ and $\f_1(x)=\f(0,x)=|x|$, for every $x\in\R$. Then, let us prove that for all $(x_1,x_2)\in\R^2$ one has $\f(x_1,x_2)=I_{[-1,1]}(x_1)+|x_2|$. This deduces from the following more general lemma, which extends observations done in the unconditional case in \cite{FGMR}.

\begin{lemma}\label{sum}
Let $n\ge1$ and $\f:\R^n\to\R\cup\{+\infty\}$ be an even convex function such that there exists two convex bodies  $K_1\subset F_1$ and $K_2\subset F_2$, where $F_1$ and $F_2$ are two complementary linear subspaces in $\R^n$, such that $\f(x_1)=\|x_1\|_{K_1}$, for all $x_1\in F_1$ and $\f(x_2)=I_{K_2}(x_2)$, for all $x_2\in F_2$. Then, $\f(x_1+x_2)=\|x_1\|_{K_1}+I_{K_2}(x_2)$, for all $x_1\in F_1$ and $x_2\in F_2$.
\end{lemma}

\begin{proof}
Let $x_1\in F_1$ and $x_2\in K_2$. From the convexity of $\f$, and using that 
\[
x_1+x_2=(1-\|x_2\|_{K_2})\frac{x_1}{1-\|x_2\|_{K_2}}+\|x_2\|_{K_2}\frac{x_2}{\|x_2\|_{K_2}},
\]
we deduce that
\[
\f(x_1+x_2)\le(1-\|x_2\|_{K_2})\f\left(\frac{x_1}{1-\|x_2\|_{K_2}}\right)+\|x_2\|_{K_2}\f\left(\frac{x_2}{\|x_2\|_{K_2}}\right)=\|x_1\|_{K_1}.
\]
On the other hand, using that $\frac{x_1}{2}=\frac{1}{2}(x_1+x_2)+\frac{1}{2}(-x_2)$, one gets
\[
\frac{\|x_1\|_{K_1}}{2}=\f\left(\frac{x_1}{2}\right)\le\frac{1}{2}\f(x_1+x_2)+\frac{1}{2}\f(-x_2)=\frac{1}{2}\f(x_1+x_2).
\]
We deduce that $\f(x_1+x_2)=\|x_1\|_{K_1}$. \\
Let $x_1\in F_1$ and $x_2\notin K_2$. Let $1<\mu<\|x_2\|_{K_2}$ and $\lambda=\mu/\|x_2\|_{K_2}\in(0,1)$. Then, $\lambda x_2\notin K_2$ and  
\[
\lambda x_2=\lambda (x_1+x_2)+(1-\lambda) \frac{-\lambda x_1}{1-\lambda}.
\]
Hence, using the convexity of $\f$, we get
\[
+\infty=\f(\lambda x_2)\le\lambda \f(x_1+x_2)+(1-\lambda) \f\left(\frac{-\lambda x_1}{1-\lambda}\right).
\]
Since $(1-\lambda) \f\left(\frac{-\lambda x_1}{1-\lambda}\right)=\lambda\|x_1\|_{K_1}<+\infty$, we deduce that $\f(x_1+x_2)=+\infty=I_{K_2}(x_2)$.
We conclude that $\f(x_1+x_2)=\|x_1\|_{K_1}+I_{K_2}(x_2)$, for all $x_1\in F_1$ and $x_2\in F_2$.
\end{proof}

{\bf B.} If $\f_2(s)=\f(se_1)=I_{[-1,1]}(s)$ and $\f_1(s)=\f(se_2)=I_{[-1,1]}(s)$, for every $s\in\R$. Let $U=\{x; \f(x)=0\}$ and $K=\dom(\f)$. From the hypothesis, one has $\f(\pm e_i)=0$, thus $\pm e_i\in U$. Since $\min\f=0$, the convexity of $\f$ implies that $U$ is convex. Thus, one has $B^{2}_{1}\subset U\subset K$. Since $\pm e_i\in\partial K$, for $i=1,2$, one deduces that there exists $u_i\in\partial K^*$ such that $\langle e_i,u_i\rangle=1$ and one has $K\subset\{x; |\langle x,u_i\rangle|\le 1, i=1,2\}$. We distinguish two cases:\\
- \underline{if $u_1=u_2$:} since $\langle e_i,u_i\rangle=1$, one has $u_1=u_2=e_1+e_2$. Thus, $K\subset\{x\in\R^2;|x_1+x_2|\le1\}:=D$. Hence, $\conv(0,e_1,e_2)\subset U\cap\R_+^2\subset K\cap\R_+^2\subset D\cap\R_+^2=\conv(0,e_1,e_2)$. Therefore, $\f_{|\R_+^2}=I_{B_1^2\cap\R_+^2}$. Using the equipartition and the fact that $\f\le I_{B_1^2}$ we conclude that $\f= I_{B_1^2}$.\\
- \underline{if $u_{1}\neq u_{2}$:} using that, for every $x\in K$, one has $\langle u_i,x\rangle\le1$ and $\f(x)\ge0$, then, for every $s>0$, and for $i=1,2$, we get
\[
\L\f(su_i)=\sup_{x}(\langle su_i,x\rangle-\f(x))=\sup_{x\in K}(s\langle u_i,x\rangle-\f(x))=s,
\]
with equality for $x=e_i$. Since $\f$ is even, we deduce that, for every $s\in\R$,
\[
\L\f(su_i)=|s|.
\]
Define $C_+=\R_+u_1+\R_+u_2$ and $C_-=\R_+u_1+\R_-u_2$. Since $e_i\in K$ and $u_i\in K^*$, one has $|\langle u_1,e_2\rangle|\le1$ and $|\langle u_2,e_1\rangle|\le1$. 
Denote by $v_i$ the unitary exterior normal of $C_+$ to the line $\R u_i$. We have $V_{\partial C_+}(\L\f)=-V_1-V_2$, where
\[
V_1=-v_2\int_{\R_+u_2}e^{-\L\f}=-v_2\int_0^{+\infty}e^{-\L\f(su_2)}ds|u_2|=-v_2\int_0^{+\infty}e^{-s}ds|u_2|=-v_2|u_2|
\]
and in the same way $V_2=-v_1|u_1|$. Hence, $V_{\partial C_+}(\L\f)=v_1|u_1|+v_2|u_2|$. It is easy to see that $V_{\partial C_+}(\L\f)\in \R_+^2$. We also have $V_{\partial\R^{2}_{+}}(\f)=-(e_1+e_2)$. Using that $\langle e_i,u_i\rangle=1$, one has 
$\langle e_1,V_1\rangle= -\langle e_1,v_2\rangle|u_2|=\langle e_2,u_2\rangle=1.$
In the same way, one also has $\langle e_2,V_2\rangle=1$. We reproduce the same argument as before with $\nabla\f(\R_+^2)$ replaced by $C_+$ and $\nabla\f(\R_+\times\R_-)$ replaced by $C_-$. Some terms are simplified because $C_+$ and $C_-$ are cones. Using again that $\f$ is even, we have
\[
P(\f)=8\int_{\R_+^{2}}e^{-\f}\int_{C_+}e^{-\L\f}+8\int_{\R_+\times\R_{-}}e^{-\f}\int_{C_-}e^{-\L\f}=8(F_1(1)+F_2(1)),
\]
where $F_1(t)=F_{\R_+^2,C_+}(t)$ and $F_2(t)=F_{\R_+\times\R_{-},C_-}(t)$.
We apply Corollary \ref{coro-deriv} and Remark~\ref{rk:asscoro} for $A=\R_+^{2}$ and $B=C_+$ and use the equipartition to get 
\begin{eqnarray}\label{ineq:B++}
\frac{d}{dt}\left(t^2F_1(t)\right)_{|t=1}\ge \langle V_{\partial \R_+^2}(\f), V_{\partial C_+}(\L\f)\rangle= \langle e_1+e_2,V_1+V_2\rangle.
\end{eqnarray} 
Applying  Corollary \ref{coro-deriv} and Remark~\ref{rk:asscoro} for $A=\R_+\times\R_{-}$ and $B=C_{-}$, we also get
\begin{eqnarray}\label{ineq:F2B}
\frac{d}{dt}\left(t^2F_2(t)\right)_{|t=1}\ge\langle V_{\partial(\R_+\times\R_-)}(\f), V_{\partial C_-}(\L\f)\rangle=\langle e_1-e_2,V_1-V_2\rangle.
\end{eqnarray} 
Adding (\ref{ineq:B++}) and (\ref{ineq:F2B}), we obtain 
\[
32=\frac{d}{dt}\left(t^2P(t\f)\right)_{|t=1}=8\frac{d}{dt}\left(t^{2}F_{1}(t)+t^{2}F_{2}(t)\right)_{|t=1}\ge 16(\langle e_1,V_1\rangle +\langle e_2,V_2\rangle)=32.
\]
Hence, we have equality in the inequalities (\ref{ineq:B++}) and (\ref{ineq:F2B}). 
From the equality case of Corollary~\ref{coro-deriv}, with $A=\R_+^2$ and $B=C_+$, we deduce that $\f$ is affine on $[-\frac{V_{\partial C_+}(\L\f)}{\mu_{\L\f}(C_+)},a]$ for every $a\in \R_+^2\cap K$. Moreover, from Proposition~\ref{prop:ineq-phi}, one has $-\frac{V_{\partial C_+}(\L\f)}{\mu_{\L\f}(C_+)}\in K\cap\R_+^2$. Since $\f$ is affine on $[\frac{V_1+V_2}{\mu_{\L\f}(C_+)},0]$ and vanishes on $B^2_1$, then $\frac{V_1+V_2}{\mu_{\L\f}(C_+)}\in U$. In the same way, we prove that $\frac{V_1-V_2}{\mu_{\L\f}(C_-)}\in U$. Hence, 
$$
\conv\left(B^2_1,\pm \frac{V_1+V_2}{\mu_{\L\f}(C_+)},\pm \frac{V_1-V_2}{\mu_{\L\f}(C_-)}\right)\subset U. 
$$ 
Then,
\begin{equation}\label{equ}
    \int_{\R^2_+}e^{-\f}\geq |U\cap\R^2_+|\geq \left|\conv\left(0,e_1,e_2,\frac{V_1+V_2}{\mu_{\L\f}(C_+)}\right)\right|=\frac{1}{2}\langle \frac{V_1+V_2}{\mu_{\L\f}(C_+)},e_1+e_2\rangle.
\end{equation}
Thus, we get $\mu_{\f}(\R^2_+)\mu_{\L\f}(C_+)\geq \frac{1}{2}\langle V_1+V_2,e_1+e_2\rangle $. Similarly, we have $\mu_{\f}(\R_+ \times\R_- )\mu_{\L\f}(C_-)\geq\frac{1}{2}\langle V_1-V_2,e_1-e_2\rangle$. Adding these two inequalities, we obtain 
$$
2=\mu_{\f}(\R^2_+)(\mu_{\L\f}(C_+)+\mu_{\L\f}(C_-))\geq \langle V_1,e_1\rangle+\langle V_2,e_2\rangle=2,
$$ 
so we get equality in (\ref{equ}). Hence, $\int_{R^2_+}e^{-\f}=|U\cap\R^2_+|$ and $\f=I_{U}=I_{K}$,
then $P(\f)=16=2P(K)$, therefore $P(K)=8$.
Thus, $K$ satisfies the equality case of Mahler's inequality in dimension 2, which implies that $K$ is a symmetric parallelogram by \cite{Re,M}.\\

{\bf C.} If $\f_2(s)=\f(se_1)=|s|$ and $\f_1(s)=\f(se_2)=|s|$, for every $s\in\R$, then $\L\f_{1}=\L\f_{2}=I_{[-1,1]}$ and $\dom(\L\f)$ is bounded. Indeed, let's prove that $\f(x)\leq\|x\|_{1}$. For all $x=(x_{1},x_{2})\in\R^{2}_{+}$,
\begin{eqnarray*}
\f(x_{1},x_{2})&=&\f\left(\frac{x_{1}}{x_{1}+x_{2}}(x_{1}+x_{2})e_{1}+\frac{x_{2}}{x_{1}+x_{2}}(x_{1}+x_{2})e_{2}\right)\\
&\leq&\frac{x_{1}}{x_{1}+x_{2}}\f((x_{1}+x_{2})e_{1})+\frac{x_{2}}{x_{1}+x_{2}}\f((x_{1}+x_{2})e_{2})
=x_{1}+x_{2}.
\end{eqnarray*}
Applying this in the other quadrants, we get that $\f(x)\leq\|x\|_{1}$ for all $x\in\R^{2}$. Hence, $\L\f(x)\geq I_{B^{2}_{\infty}}(x)$ and $\dom(\L\f)\subset B^{2}_{\infty}$ is bounded. Thus, there exists a linear invertible map $T$ such that $\psi=(\L\f)\circ T$ is strongly equipartitioned  and $P(\psi)=P(\L\f)=16$. Then, for all $i=1,2$, $P(\psi_{i})=4$ and $\psi_{i}(x)=I_{[-1,1]}(x)$ or $|x|$. Since $\dom(\psi)$ is bounded, then $\psi_{i}=I_{[-1,1]}$. From case B, one concludes that $\psi=I_{K}$, where $K$ is a symmetric  parallelogram hence $\L\f=I_{L}$, where $L=T^{-1}(K)$ is a symmetric parallelogram.
\end{proof}
\section{The inequality in dimension $n$}
\begin{thm}\label{thm:dimn}
Let $\f:\R^{n}\to\R\cup\{+\infty\}$ be even convex such that $\f$ and $\f_{i}$ are equipartitioned for all $1\leq i\leq n$ and  $P(\f_{i})\geq 4^{n-1}$
then $P(\f)\geq 4^{n}$.
\end{thm}
\begin{proof}
We can assume that $\f(0)=0$. From Lemma \ref{regul}, we reduce to the case where $\dom(\f)=\dom(\L\f)=\R^{n}$, $\f$ is differentiable and strictly convex on $\R^{n}$ and $\R^{n}=\cup_{\varepsilon\in\{-1,1\}^{n}}\nabla\f(\R^{n}_{\varepsilon})$. Using the equipartition, we have 
$$
P(\f)=\int_{\R^{n}}e^{-\f}\int_{\R^{n}}e^{-\L\f}=2^{n}\int_{\R^{n}_{\varepsilon}}e^{-\f}\int_{\R^{n}}e^{-\L\f}=2^{n}\sum_{\varepsilon\in\{-1,1\}^{n}}F_{\varepsilon}(1),
$$
where $F_{\varepsilon}(t)=F_{\R^{n}_{\varepsilon},\nabla\f(\R^{n}_{\varepsilon})}(t)$, for every $\varepsilon\in\{-1,1\}^{n}$. Using 4) of Lemma \ref{cor:deriv} and Remark \ref{rk:asscoro}, for $A=\R^{n}_{\varepsilon}$ and $B=\nabla\f(\R^{n}_{\varepsilon})$, we get that, for every $\varepsilon\in\{-1,1\}^{n}$, one has 
$$
F_{\varepsilon}'(1)=nF_{\varepsilon}(1)-\int_{\R^{n}_{\varepsilon}}e^{-\f}\int_{\nabla\f(\R^{n}_{\varepsilon})}\L\f e^{-\L\f}-\int_{\R^{n}_{\varepsilon}}\f e^{-\f}\int_{\nabla\f(\R^{n}_{\varepsilon})}e^{-\L\f}.
$$
Thus, using the equipartition, we have
\begin{eqnarray*}
\frac{d}{dt}\left(t^{n}F_{\varepsilon}(t)\right)_{|t=1}&=&nF_{\varepsilon}(1)+F_{\varepsilon}'(1)\\
&=&2nF_{\varepsilon}(1)-\int_{\R^{n}_{\varepsilon}}e^{-\f}\int_{\nabla\f(\R^{n}_{\varepsilon})}\L\f e^{-\L\f}-\int_{\R^{n}_{\varepsilon}}\f e^{-\f}\int_{\nabla\f(\R^{n}_{\varepsilon})}e^{-\L\f}\\
&=&2nF_{\varepsilon}(1)-\frac{1}{2^{n}}\left(\int_{\R^{n}}e^{-\f}\int_{\nabla\f(\R^{n}_{\varepsilon})}\L\f e^{-\L\f}-\int_{\R^{n}}\f e^{-\f}\int_{\nabla\f(\R^{n}_{\varepsilon})}e^{-\L\f}\right).
\end{eqnarray*}
Summing these terms and using again the equipartition, we get
\begin{eqnarray*}
\sum_{\varepsilon}\frac{d}{dt}\left(t^{n}F_{\varepsilon}(t)\right)_{|t=1}&=&2n\sum_{\varepsilon}F_{\varepsilon}(1)-\frac{1}{2^{n}}\sum_{\varepsilon}\left(\int_{\R^{n}}e^{-\f}\int_{\nabla\f(\R^{n}_{\varepsilon})}\L\f e^{-\L\f}-\int_{\R^{n}}\f e^{-\f}\int_{\nabla\f(\R^{n}_{\varepsilon})}e^{-\L\f}\right)\\
&=&2n\frac{P(\f)}{2^{n}}-\frac{1}{2^{n}}\left(\int_{\R^{n}}e^{-\f}\int_{\R^{n}}\L\f e^{-\L\f}-\int_{\R^{n}}\f e^{-\f}\int_{\R^{n}}e^{-\L\f}\right).
\end{eqnarray*}
Now, applying Lemma \ref{cor:deriv} and Remark \ref{rk:asscoro}, one has
\begin{eqnarray*}
\frac{d}{dt}\left(t^{n}P(t\f)\right)_{|t=1}&=&nP(\f)+\left(nP(\f)-\int_{\R^{n}}e^{-\f}\int_{\R^{n}}\L\f e^{-\L\f}-\int_{\R^{n}}\f e^{-\f}\int_{\R^{n}}e^{-\L\f}\right)\\
&=&2^{n}\sum_{\varepsilon}\frac{d}{dt}\left(t^{n}F_{\varepsilon}(t)\right)_{|t=1}.
\end{eqnarray*}
We apply Corollary \ref{coro-deriv} and Remark \ref{rk:asscoro}  for $A=\R^{n}_{\varepsilon}$ and $B=\nabla\f(\R^{n}_{\varepsilon})$ and use the equipartition to get
\[
\frac{d}{dt}\left(t^{n}P(t\f)\right)_{|t=1}\geq2^{n}\sum_{\varepsilon}\left(\left\langle V_{\partial\R^{n}_{\varepsilon}}(\f),V_{\partial\nabla\f(\R^{n}_{\varepsilon})}(\L\f)\right\rangle+\frac{\mu_{\f}(\R^{n})}{2^{n}}Q_{\partial\nabla\f(\R^{n}_{\varepsilon})}(\L\f)\right).
\]
Notice that
$$
\sum_{\varepsilon}Q_{\partial\nabla\f(\R^{n}_{\varepsilon})}(\L\f)=\sum_{\varepsilon}\sum_{i=1}^{n}Q_{\nabla\f(\R^{n}_{\varepsilon}\cap e_{i}^{\perp})}(\L\f)=\sum_{i=1}^{n}\sum_{\varepsilon}\int_{\nabla\f(\R^{n}_{\varepsilon}\cap e_{i}^{\perp})}\langle y,n_{\nabla\f(\R^{n}_{\varepsilon}\cap e_{i}^{\perp})}(y)\rangle e^{-\L\f(y)}dy,
$$
where in $\nabla\f(\R^{n}_{\varepsilon}\cap e_{i}^{\perp})$, the normal is chosen exterior to $\nabla\f(\R^{n}_{\varepsilon})$. Since in each hyperplane $e_i^\perp$, each cone $\R_\varepsilon^n\cap e_i^\bot$ appears twice with two opposite orientations, the sum of these two terms vanishes. Hence, the whole sum vanishes. Using that, in each $e_i^\bot$, the function $\f_i$ is equipartitioned and Lemma~\ref{lem:V-R+}, it follows that
\begin{eqnarray*}
\frac{d}{dt}\left(t^{n}P(t\f)\right)_{|t=1}&\geq&2^{n}\sum_{\varepsilon}\sum_{1\le i,j\le n}\int_{\R_\varepsilon^n\cap e_i^\bot}e^{-\f_i}\langle -\varepsilon_ie_i,V_{\nabla\f(\R_\varepsilon^n\cap e_j^\bot)}(\L\f)\rangle\\
&\geq&2\sum_{1\le i,j\le n}\sum_{\varepsilon}\int_{e_i^\bot}e^{-\f_i}\langle -\varepsilon_ie_i,V_{\nabla\f(\R_\varepsilon^n\cap e_j^\bot)}(\L\f)\rangle.
\end{eqnarray*}
Noticing again that, in each hyperplane $e_j^\perp$, each cone $\R_\varepsilon^n\cap e_j^\bot$ appears twice, with two opposite orientations, one has, for every fixed $i\neq j$, 
\[
\sum_{\varepsilon}\varepsilon_iV_{\nabla\f(\R_\varepsilon^n\cap e_j^\bot)}=0.
\]
Thus, 
\[
\frac{d}{dt}\left(t^{n}P(t\f)\right)_{|t=1}\geq 2\sum_{i=1}^n\int_{e_i^\bot}e^{-\f_i}\sum_{\varepsilon}\langle -\varepsilon_ie_i,V_{\nabla\f(\R_\varepsilon^n\cap e_i^\bot)}(\L\f)\rangle.
\]
Using Lemma~\ref{lem:V-R+}, we get 
\[
\frac{d}{dt}\left(t^{n}P(t\f)\right)_{|t=1}\geq 2\sum_{i=1}^{n}\int_{e_{i}^{\perp}}e^{-\f_{i}}\sum_{\varepsilon}\int_{\R^{n}_{\varepsilon}\cap e_{i}^{\perp}}e^{-\L(\f_{i})}
=2\sum_{i=1}^{n}\int_{e_{i}^{\perp}}e^{-\f_{i}}\left(2\int_{e_{i}^{\perp}}e^{-\L(\f_{i})}\right)=4\sum_{i=1}^{n}P(\f_i).
\]
Since $P(\f_{i})\geq 4^{n-1}$, we get
$$
\frac{d}{dt}\left(t^{n}P(t\f)\right)_{|t=1}\ge 4^{n}n.
$$
Applying this to $s\f$, we deduce that, for all $t>0$, 
$$
\frac{d}{dt}(t^{n}P(t\f))\geq 4^{n}nt^{n-1}.
$$
Integrating this inequality, we conclude that, for every $0<\e<t$, one has
\begin{eqnarray*}
t^nP(t\f)= \e^nP(\e\f)+\int_\e^t\frac{d}{ds}\left(s^nP(s\f)\right)ds\ge 4^{n}n\int_\e^t s^{n-1}ds=4^{n}(t^{n}-\e^{n}).
\end{eqnarray*}
Letting $\e$ tends to $0$, we conclude that $t^{n}P(t\f)\ge 4^{n}t^n$ and thus $P(\f)\ge4^n$.\\
\end{proof}

\noindent {\bf Acknowledgments:}  The authors  are grateful to the anonymous referees for a careful reading of the paper and constructive comments and corrections and they also thank Olivier Gu\'edon and St\'ephane Sabourau for their remarks and pertinent questions on this work.

\noindent
{\footnotesize\sc Matthieu Fradelizi:}
  {\footnotesize Univ Gustave Eiffel, Univ Paris Est Creteil, CNRS, LAMA UMR8050 F-77447 Marne-la-Vallée, France. }\\[-1.3mm]
  {\footnotesize e-mail: {\tt matthieu.fradelizi@univ-eiffel.fr \tt }}\\

\noindent
{\footnotesize\sc Elie Nakhle:}
  {\footnotesize Univ Paris Est Creteil, Univ Gustave Eiffel, CNRS, LAMA UMR8050, F-94010 Creteil, France. }\\[-1.3mm]
  {\footnotesize e-mail: {\tt elie.nakhle@u-pec.fr \tt }}\\

\end{document}